\documentclass[journal]{IEEEtran}

\usepackage{hyperref}
\usepackage{amssymb}
\usepackage[utf8]{inputenc} 
\usepackage[normalem]{ulem}
\usepackage[T1]{fontenc}
\usepackage{url}
\usepackage{ifthen}
\usepackage{multirow}
\usepackage{amsmath,amsthm}
\usepackage{graphicx,color,epsfig,rotating,subfigure}       
\usepackage{algorithm,algcompatible,xcolor,cite}
\usepackage{graphicx,color,epsfig,rotating,subfigure}
\usepackage{amsfonts,amsmath,amssymb, bm}
\usepackage{algorithm}
\usepackage{subfigure}
\usepackage{algpseudocode}

\long\def\comment#1{}

\newfont{\bbb}{msbm10 scaled 700}

\newfont{\bb}{msbm10 scaled 1100}

\newcommand{\RR}{\mbox{\bb R}}

\renewcommand{\arg}{{\hbox{arg}}}

\newcommand{\trasp}{{\sf T}}

\usepackage{graphicx}
\usepackage{tabularx,booktabs,hhline}
\usepackage{color, colortbl}
\definecolor{LightCyan}{rgb}{0.88,1,1}
\definecolor{lightgray}{gray}{0.95}
\usepackage{multirow}
\usepackage{soul} 
\usepackage{amsthm}

\usepackage{tikz}
\usetikzlibrary{arrows.meta, positioning}
\usepackage{booktabs}

\newcommand{\qhats}{\hat{\mathcal{Q}}\mathbf{S}}

\newcommand{\clambda}{\mathcal{C}\hat{\mathbf{\Lambda}}_{\Omega}}

\newcommand{\xcheck}{\check{\mathbf{X}}_{1}}

\newcommand{\xbar}{\bar{\mathbf{X}}_{1}}

\newtheorem{theorem}{Theorem}
\newtheorem{definition}{Definition}
\newtheorem{lemma}{Lemma}

\newcommand{\argmin}{\operatornamewithlimits{argmin}}

\renewcommand{\subsubsection}[1]{\noindent {\bf #1. }}

\begin{document}

\title{Structured Tensor Approximation from Lateral Slice Sampling via Basis and Manifold Priors} 

\author{\vspace{-0.1in}
\IEEEauthorblockN{Jeongmin Chae$^{\dagger}$, Usama Saleem$^{*}$, Selin Bac$^{\dagger\dagger}$,
Shaama Mallikarjun Sharada$^{*}$ and Urbashi Mitra$^{\dagger}$ } 
\thanks{$^{\dagger}$ J. Chae and U.Mitra are with the Department of Electrical and Computer Engineering, University of Southern California, Los Angeles, USA. E-mail: \{chaej, ubli\}@usc.edu.}
\thanks{$^{\dagger\dagger}$ S. Bac is with the Department of Chemical Engineering, University of California, Santa Barbara, USA.}
\thanks{
$^{*}$ U.Saleem and S. Mallikarjun Sharada are with the Department of Chemical Engineering and Materials Science, University of Southern California, Los Angeles, USA. E-mail: \{usaleem, ssharada\}@usc.edu.}}

\twocolumn
\maketitle
\vspace{-0.4in}
\begin{abstract}
    In this work, we consider a structured tensor approximation problem, where only a limited number of lateral slices are observed. The proposed algorithm , called Basis and Manifold prior Tensor Approximation (BMTA), exploits both global and local structures of the evolution of a global tensor. Specifically, BMTA integrates two signal models: (i) a quasi-basis model that captures smooth global variations along a physical trajectory, and (ii) a manifold-guided interpolation model that characterizes local relationships among tensor slices. A low-rank Tucker reconstruction framework is incorporated to efficiently capture the priors, resulting in coefficients for basis function estimation and a tensor optimization. In addition, we provide a theoretical analysis which establishes a non-asymptotic reconstruction error bound that characterizes the effects of sampling complexity, optimization convergence, and model mismatch. Numerical experiments are performed on both synthetic and real-world datasets, including quantum chemistry and spatiotemporal sensing applications. 
\end{abstract}

\vspace{-0.2in}
\section{Introduction}

Multidimensional data represented as tensors arise in numerous scientific and engineering applications, including medical imaging \cite{ahmadi2021cross}, sensor networks \cite{sekar2022compressed}, neuroscience \cite{li2018tucker}, and quantum chemistry \cite{chae2025quasi,chae2024matrix,bac2025incorporating}. Recovering tensors from incomplete observations is therefore a fundamental problem. Existing tensor completion and approximation methods typically exploit low-rank structure under random entry, fiber, or slice sampling assumptions \cite{kolda2009tensor}. However, such assumptions are often violated in practice, where data acquisition is constrained by physical processes and follows structured sampling patterns \cite{chae2026sketched,song2016sublinear,zhang2020high}.

A motivating example arises in quantum chemistry, where nuclear Hessian matrices along a reaction pathway are required for accurate prediction of reaction rates using variational transition state theory with multidimensional tunneling (VTST-MT) \cite{bac2022matrix,fernandez2007variational}. Since Hessian evaluations are computationally expensive, only a small number can be computed, resulting in lateral slice sampling along the reaction coordinate rather than random observations \cite{corchado1998interpolated}. This motivates reconstructing the full tensor from a limited number of structured slice observations.

Empirical evidence suggests that Hessian tensors exhibit two complementary structures: frontal slices evolve smoothly along the reaction coordinate and are well represented by structured basis functions, such as Legendre polynomial or discrete cosine transform (DCT) bases \cite{corchado1998interpolated}, while neighboring slices exhibit strong local similarity, indicating an underlying low-dimensional manifold. These observations motivate incorporating both global and local structural priors into tensor approximation. 

Motivated by these observations, we propose the {\em Basis and Manifold prior Tensor Approximation (BMTA)} framework, which integrates (i) a structured basis representation (i.e., Legendre polynomial and Discrete Cosine Transform (DCT) basis) describing the global evolution of frontal slices and (ii) a manifold-guided interpolation model capturing local geometric relationships between neighboring slices. Unlike conventional tensor completion methods, BMTA is designed for structured lateral slice sampling and combines closed-form basis coefficient estimation with constrained low-rank Tucker optimization. In the numerical section, we show BMTA accommodates various structured basis representations.

From a theoretical perspective, we establish a non-asymptotic reconstruction error bound that characterizes the effects of optimization, sampling complexity, and model mismatch. The bound decomposes the reconstruction error into a contracting optimization term and an irreducible approximation error induced by the structural priors. Extensive experiments on synthetic and real-world datasets demonstrate that BMTA consistently outperforms existing approaches under limited observations while remaining robust to noise and model mismatch.

\begin{itemize}
\item {\bf Structured tensor approximation.}
We propose BMTA, which integrates a structured basis representation and manifold-guided priors for tensor reconstruction from lateral slice sampling.

\item {\bf Theoretical guarantees.}
We derive a non-asymptotic reconstruction error bound characterizing sampling complexity, optimization convergence, and model mismatch.

\item {\bf Experimental validation.}
Experiments on synthetic and real-world datasets demonstrate significant improvements over existing methods, particularly under limited observations.
\end{itemize}
\vspace{-0.2in}
\section{Preliminaries}\label{sec:prelim}
We employ the following notation. Tensors and matrices are denoted by $\mathcal{X}$ and $\mathbf{X}$, respectively. We consider third-order tensors motivated by the applications in \cite{quiton2022toward,chae2024matrix,chae2025quasi}; the proposed framework naturally extends to higher-order tensors. Slices are second-order tensors obtained by fixing all modes except two. For a third-order tensor $\mathcal{X}$, the slices $\mathcal{X}(:,:,k)$, $\mathcal{X}(:,j,:)$, and $\mathcal{X}(i,:,:)$ are referred to as the \emph{frontal}, \emph{lateral}, and \emph{horizontal} slices, respectively.

The Frobenius norm of a tensor is
$
\|\mathcal{X}\|_F^2
=
\langle\mathcal{X},\mathcal{X}\rangle
=
\sum_{i_1,i_2,i_3}
\mathcal{X}(i_1,i_2,i_3)^2,$
and the spectral norm is
$
\|\mathcal{X}\|
=
\sup_{\|{\bf u}_k\|_2\le1}
\big|
\langle
\mathcal{X},
{\bf u}_1\circ{\bf u}_2\circ{\bf u}_3
\rangle
\big|,$
where ${\bf u}_1\circ{\bf u}_2\circ{\bf u}_3$ denotes a rank-one tensor. Let $GL(r)$ denote the set of invertible $r\times r$ real matrices. For a matrix $\mathbf{A}$, $\|\mathbf{A}\|$ and $\|\mathbf{A}\|_F$ denote its spectral and Frobenius norms, respectively. Finally, $f(n)=O(g(n))$ if there exists a constant $c>0$ such that
$
|f(n)|\le c|g(n)|,  \forall n\in\mathbb{N}.$
We next introduce key definitions to be used throughout the paper.
\begin{figure}[t]
    \centering
    \includegraphics[width=0.9\linewidth]{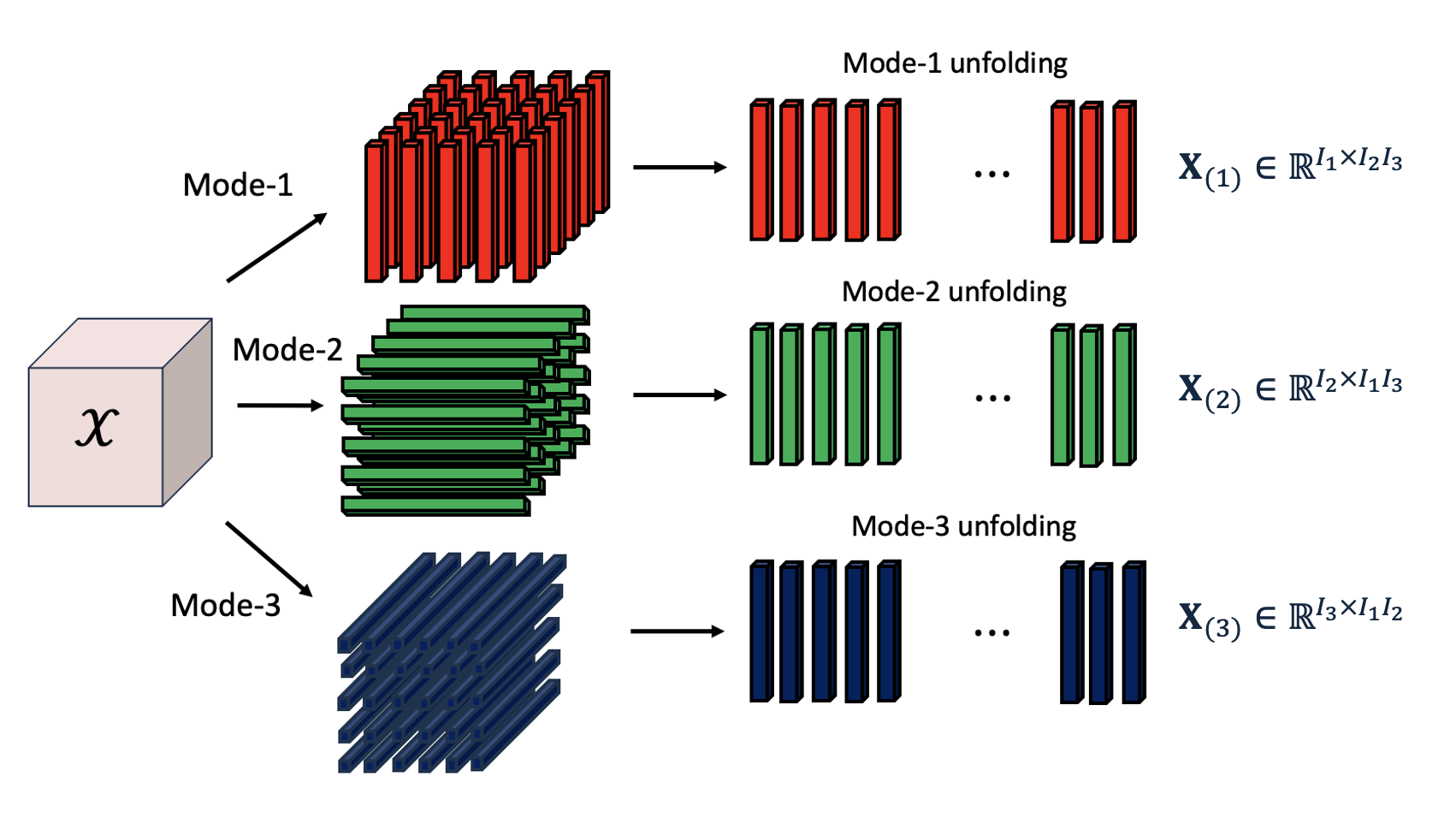}
    \vspace{-0.15in}
    \caption{Illustration of mode-(n) tensor unfolding and folding for a third-order tensor.}
    \label{fig:def1}
      \vspace{-0.25in}
\end{figure}

\begin{definition}{\bf (Mode-$n$ Tensor unfolding and folding)} Given a tensor $\mathcal{X} \in \RR^{I_{1} \times \dots \times I_{N}}$, we denote the mode-$n$ unfolding of $\mathcal{X}$ as $\mathbf{X}_{(n)}\in \mathbb{R}^{I_n \times \Pi_{1\leq t\leq N, t\neq n} I_t}$ , whose $(i_{n},j)$-th element in the lexicographical order is mapped from the $(i_{1},i_{2},\dots,i_{N})$-th entry of tensor $\mathcal{X}$, where $j=1+\sum_{1\leq l \leq N, l \neq n} (i_{l}-1)J_{l}\; \mbox{with}\; J_{l} = \Pi_{1\leq t \leq l-1, t\neq n} I_{t}$. The inverse operator to express the tensor $\mathcal{X}$, is defined as $\mathcal{X}=\mbox{fold}_{n}\left(\mathbf{X}_{(n)}\right)$ (See Figure~\ref{fig:def1}).
\end{definition}
\vspace{-0.1in}
\begin{definition}
    {\bf (Multilinear rank)} {Let $\text{rank}(\mathbf{X})$ denote the standard matrix rank of $\mathbf{X}$, i.e., the dimension of the column space of a matrix.} Given a N order tensor $\mathcal{X} \in \RR^{I_{1} \times \dots \times I_{N}}$, the multilinear rank of $\mathcal{X}$ is defined as the tuple $\text{rank}_{ML}(\mathcal{X})=(r_1,r_2,\dots,r_N)$, where $r_n = \text{rank}(\mathbf{X}_{(n)})$ for all $n\in[N]$. Equivalently, $r_n$ is the dimension of the vector space spanned by the mode-$n$ fibers of $\mathcal{X}$. 
\end{definition}
{For example, in Figure~\ref{fig:def1}, multilinear rank of $\mathcal{X}$ is defined as $(r_1,r_2,r_3)$, where $r_1$ = $\text{rank}\left(\mathbf{X}_{(1)}\right)$, $r_2$ = $\text{rank}\left(\mathbf{X}_{(2)}\right)$ and $r_3$ = $\text{rank}\left(\mathbf{X}_{(1)}\right)$.}

\begin{figure}[t]
    \centering
    \includegraphics[width=0.8\linewidth]{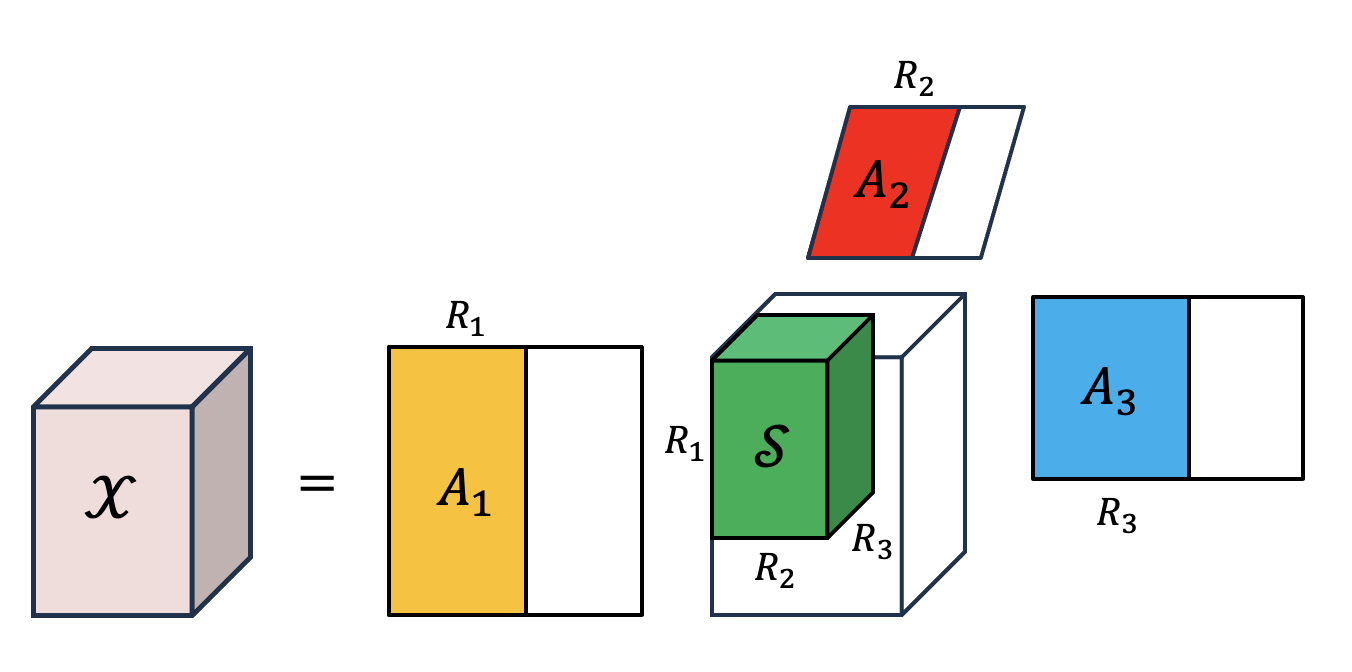}
    \vspace{-0.15in}
    \caption{Illustration of Tucker decomposition for a third-order tensor.}
    \label{fig:def3}
       \vspace{-0.25in}
\end{figure}

\begin{definition}
    {\bf (Tucker decomposition \cite{kolda2009tensor})} Let $\mathcal{X} \in \RR^{I_{1} \times I_{2} \times \dots \times I_{N}}$ be a given data tensor, then the Tucker decomposition model is given by 
    $
    \mathcal{X} = \mathcal{S} \times_{1} {\mathbf{A}_{1}} \times_{2} \mathbf{A}_{2} \cdots \times_{N} \mathbf{A}_{N},$
 where $\mathcal{S} \in \RR^{R_{1} \times R_{2} \times \dots \times R_{N}}$ is called the core tensor, and the matrices $\mathbf{A}_{n} \in \RR^{I_{n} \times R_{n}}$, $R_{n} \leq I_{n}$, $n =1,2, \dots, N$ are called \emph{factor matrices}; they are multiplied by the core tensor $\mathcal{S}$ along each mode $n$. In short, we have $\mathcal{X}=(\mathbf{A}_{1},\mathbf{A}_{2},\dots,\mathbf{A}_{N})\cdot \mathcal{S}$ (See Figure~\ref{fig:def3}).
 \end{definition}
We define the largest and smallest nonzero singular values across all mode-$k$ unfoldings as
$
\sigma_{\max}(\mathcal{X})
=
\max_k \sigma_{\max}\big(\mathbf{X}_{(k)}\big),
\qquad
\sigma_{\min}(\mathcal{X})
=
\min_k \sigma_{\min}\big(\mathbf{X}_{(k)}\big).$
Note that, in general, $\|\mathcal{X}\|\neq\sigma_{\max}(\mathcal{X})$.

For a tensor $\mathcal{X}\in\mathbb{R}^{I_1\times\cdots\times I_N}$ and a matrix $\mathbf{B}\in\mathbb{R}^{J\times I_n}$, the mode-$n$ product is denoted by
$\mathcal{Y}=\mathcal{X}\times_n\mathbf{B}$, whose mode-$n$ unfolding satisfies
$\mathbf{Y}_{(n)}
=
\mathbf{B}\mathbf{X}_{(n)}.$
Furthermore, if
$
\mathcal{X}
=
\mathcal{S}
\times_1\mathbf{A}_1
\times_2\mathbf{A}_2
\cdots
\times_N\mathbf{A}_N,$
then
\begin{align}
\operatorname{vec}(\mathcal{X})
&=
\left(\bigotimes_{n=1}^{N}\mathbf{A}_n\right)
\operatorname{vec}(\mathcal{S}),
\label{eq:Kronecker}\\
\mathbf{X}_{(n)}
&=
\mathbf{A}_n
\mathbf{S}_{(n)}
\left(
\bigotimes_{i\neq n}\mathbf{A}_i
\right)^{\top},
\label{eq:Kronecker_unfolding}
\end{align}
where
$\bigotimes_{n=1}^{N}\mathbf{A}_n
=
\mathbf{A}_N\otimes\cdots\otimes\mathbf{A}_1$.
\vspace{-0.1in}
\section{Problem Formulation and Optimization}
\subsection{System model}
Let $\mathcal{H} \in \RR^{m \times n \times m}$ be the ground truth tensor of $n$ different  matrices of size $m \times m$.  {For our chemical engineering motivating application, } the lateral slices $\mathcal{H}(:,i_2,:)$, $1\leq i_2\leq n$, would be the Hessian matrices generated from the corresponding chemical reaction systems at time $i_2 \in {[\mathbf{s}]}$, where $\mathbf{s}=[1,\dots,s]$ denotes the ordered sequence of reaction times.  We assume the ground truth tensor $\mathcal{H}$ admits the following Tucker decomposition, for $1\leq i_1,i_3 \leq m$ and $1\leq i_2 \leq n$,
\begin{align}\label{eq:Tucker}
    &\mathcal{H}(i_1,i_2,i_3)\\ &= \sum_{j_{1}=1}^{r_1}\sum_{j_{2}=1}^{r_2}\sum_{j_{3}=1}^{r_1}\; \mathbf{X}_{1}(i_1,j_1)\;\mathbf{X}_{2}(i_2,j_2)\mathbf{X}_{1}(i_3,j_3)\;\mathcal{G}(j_1,j_2,j_3),\nonumber
\end{align} namely, $\mathcal{H}=(\mathbf{X}_{1},\mathbf{X}_{2},\mathbf{X}_{1})\cdot \mathcal{G}$,
where $\mathcal{G} \in \RR^{r_1\times r_2 \times r_1}$ is the core tensor of multilinear rank $\mathbf{r}=(r_1,r_2,r_1),$ $r_{k} = \mbox{rank}(\mathbf{H}_{(k)})$ for $k=1,2$. $\mathbf{X}_{1} \in \RR^{m \times r_1}$, and $\mathbf{X}_{2} \in \RR^{n \times r_2}$ are factor matrices of each mode, where $1 \leq r_1 \leq m$ and $1 \leq r_2 \leq n$. The factor matrices of mode-1 and mode-3 are the same ($\mathbf{X}_{1}$) due to the symmetry of a Hessian matrix. We fix the ground truth factors $\mathbf{X}_1$ and $\mathbf{X}_2$ as orthonormal matrices consisting of left singular vectors in each mode. Let $\mathbf{G}_{(k)}$ be  the mode-$k$ matricization of $\mathcal{G}$.  The core tensor $\mathcal{G}$ is related to the singular values in each mode as
\begin{align}
    \mathbf{\Sigma}^2_k = \mathbf{G}_{(k)}\mathbf{G}_{(k)}^{\trasp}, \;\; k=1,2, 
\end{align} where $\mathbf{\Sigma}_k=\text{diag}\left[\sigma_1\left(\mathbf{H}_{(k)})\right),\dots,\sigma_{r_k}\left(\mathbf{H}_{(k)})\right)\right]$ is a diagonal matrix where the diagonal elements are composed of the nonzero singular values of $\mathbf{H}_{(k)}$ and $r_k=\text{rank}\left(\mathbf{H}_{(k)}\right)$ for $k=1,2$.

\begin{figure}[t]
    \centering
    \includegraphics[width=0.8\linewidth]{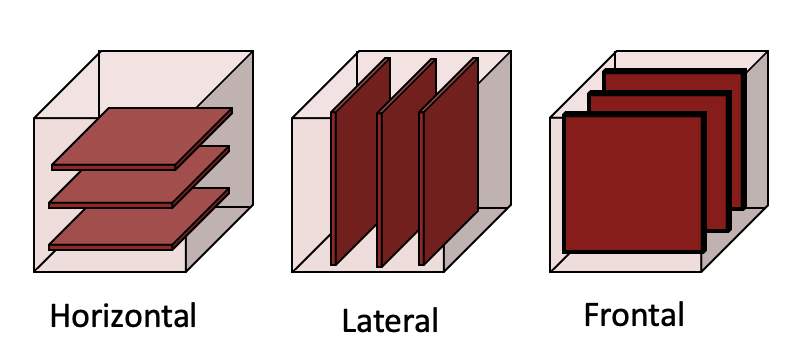}
    \vspace{-0.15in}
    \caption{Illustration of various slice sampling models.}
    \label{fig:lateral}
        \vspace{-0.25in}
\end{figure}

\subsubsection{Lateral slice sampling}  Our goal is to reconstruct the true tensor $\mathcal{H}$ only given $d$ sampled lateral slices of $\mathcal{H}$. We assume that the $d$ lateral slices are sampled uniformly at random. 
Let $\mathbf{\Psi}$ be the sampling matrix for lateral slices, which samples $d$ lateral slices from a total of $n$. Let $\Omega = \{i_{1},\dots, i_{d}\} \subset [n]$ denote the set of sampled
lateral slice indices. Clearly $|\Omega| = d$. Then, $\mathbf{\Psi} \in \{0,1\}^{d \times n}$ is $\mathbf{\Psi} \doteq \mathbf{I}_{\Omega}$,
 where $\mathbf{I}$ is the identity matrix of dimension $n$ and the notation
$\mathbf{I}_{\Omega}$ means that we consider the sub-matrix of $\mathbf{
I}$ formed by its
rows indexed by entries in the set $\Omega$. Therefore, a tensor of $d$ sampled lateral slices $\mathcal{C} \in \RR^{m \times d \times m}$ can be expressed as $\mathcal{C}=\mbox{fold}_{2}\left( \mathbf{\Psi}\mathbf{H}_{(2)}\right).$ Examples of various slicing strategies are provided in Figure~\ref{fig:lateral}.

\subsubsection{Low-dimensional Basis Representation}
Many physical systems exhibit low-dimensional structure along an underlying trajectory, such as reaction coordinates, time, or spatial parameters \cite{coifman2008diffusion,bac2022matrix}.  Motivated by this observation, we assume that the frontal slices of $\mathcal{H}$ admit a low-dimensional basis representation along the second mode,
\begin{align}
\label{eq:polynomial}
\mathcal{H}
=
\mathcal{Q}
\times_{1}\mathbf{I}
\times_{2}\mathbf{S}
\times_{3}\mathbf{I}
+
\mathcal{E}_{QS},
\end{align}
where $\mathbf{S}\in\mathbb{R}^{n\times l}$ is a known basis matrix, $\mathcal{Q}\in\mathbb{R}^{m\times l\times m}$ is the coefficient tensor, and $\mathcal{E}_{QS}\in\mathbb{R}^{m\times n\times m}$ models the representation error. Here, $l$ denotes the basis dimension and $\mathbf{I}\in\mathbb{R}^{m\times m}$ is identity matrix.

Model \eqref{eq:polynomial} is a perturbed Tucker representation in which $\mathcal{Q}$ acts as the core tensor and $\mathbf{S}$ is the factor matrix associated with the trajectory mode. The choice of $\mathbf{S}$ is application dependent and may correspond to Legendre polynomial bases for smoothly varying trajectories \cite{chae2024matrix,platte2005polynomials}, DCT/Fourier bases for temporal signals \cite{strang1999discrete,chen2024double}, or other physics-informed basis representations.

\subsubsection{Manifold-guided Slice Interpolation}  
While the low-dimensional basis representation in \eqref{eq:polynomial} captures the global evolution of tensor slices, the exact relationship between slices along the underlying trajectory is generally unknown, and $\mathbf{S}$ may not fully represent local variations in the trajectory. Nevertheless, in many physical systems, neighboring states along the trajectory exhibit strong correlations, suggesting an underlying smooth manifold structure.
Motivated by
this, we introduce an interpolative structure that leverages the
pairwise distances between slices to reconstruct the full tensor from a small set of sampled slices leveraged by the reaction coordinate $\mathbf{s}$.

Given the set of $d$ sampled lateral slice indices $\Omega=\{i_{1},\dots,i_{d}\} \subset [n]$ and the tensor of sampled lateral slices $\mathcal{C} \in \RR^{m \times d \times m}$, we adopt the following model:
\vspace{-0.03in}
\begin{align}\label{eq:interpolatory}
    \mathcal{H}=\mathcal{C} \times_{1} \mathbf{I} \times_{2} \mathbf{\Lambda}_{\Omega}\times \mathbf{I} + \mathcal{E}_{\mathcal{C}},
\end{align}
where $\mathbf{\Lambda}_{\Omega}\in \left[0,1\right]^{n\times d}$ is a structured interpolation coefficient matrix with respect to each index in $\Omega$. 
$\mathbf{\Lambda}_{\Omega}$ models the change as the reaction is traversed.
$\mathcal{E}_{\mathcal{C}}$ is an unknown perturbation tensor with respect to $\mathcal{C}$. 
The interpolation model is motivated by the assumption that the tensor slices evolve on a smooth underlying manifold parameterized by $\mathbf{s}$. Since the exact mapping from the trajectory coordinate to the tensor space is unknown, we approximate this relationship by constructing $\mathbf{\Lambda}_{\Omega}$ from pairwise distances between trajectory points, as described in the following section.
\vspace{-0.2in}
\subsection{Optimization} 
The proposed algorithm consists of three stages. First, the coefficient tensor $\mathcal{Q}$ is estimated via the basis model in \eqref{eq:Qhat}. Second, the interpolation matrix $\mathbf{\Lambda}_{\Omega}$ is constructed according to \eqref{eq:weight}. Finally, the Tucker factors are optimized by minimizing \eqref{eq:optimization}, jointly enforcing the basis representation and manifold-guided interpolation constraints.

\subsubsection{1. Quasi-structure  estimation} We first estimate the structured basis coefficients tensor $\mathcal{Q}$ using $\mathcal{C} = \mbox{fold}_{2}\left(\mathbf{\Psi}\mathbf{H}_{(2)}\right)$. We argue that as long as enough independent lateral slices are sampled (Theorem~\ref{thm:main}), the following optimization
gives us a good estimate for $\mathcal{Q}$, 
\begin{align}\label{eq:Qhat}
    \mathcal{\hat{Q}} = \argmin_{\mathcal{\bar{Q}}}\;\left\Vert \mathbf{C}_{(2)}-\mathbf{\Psi}\left(\mathcal{\bar{Q}} \times_{2} \mathbf{S}\right)_{(2)}\right\Vert_{F}^{2},
\end{align}
where $\mathbf{C}_{(2)}$ is the mode-2 unfolding of $\mathcal{C}$. This is a standard tensor regression problem {\cite{ahmed2020tensor,cai2019nonconvex}} that admits a closed
form solution, where 
$\mathcal{Q}\doteq\mbox{fold}_{(2)}\left(\mathbf{\Psi}^{\dagger}\mathbf{C}_{(2)}\right) \times_{2} \mathbf{S}^{\dagger}$, where $\dagger$ denotes Moore–Penrose inverse.  This step exploits the fact that lateral slide side information, $\mathbf{S}$, is known.

\subsubsection{2. Interpolative tensor estimation}
{We next determine the interpolation weights $\mathbf{\Lambda}_{\Omega}$ based on the pairwise distances calculated using $\mathbf{s}=[s_{1}, s_{2}, \dots, s_{n}]$ between slices along the reaction coordinate.} $s_j$ denotes the reaction coordinate value associated with the 
$j$th lateral slice. We define the distance between slices $j$ and $i$ as$\left\Vert s_j-s_i\right\Vert^2$. We then construct interpolation weights using a Gaussian (RBF) kernel \cite{jayasumana2015kernel}. 

Specifically, for each unsampled index $j \in [n] \setminus \Omega$ and sampled index $i \in \Omega$, we compute kernel weights as
\begin{align}\label{eq:weight}
\tilde w_{j,i}
=
\exp\!\left(
-\frac{\left\Vert s_j-s_i\right\Vert^2}{2\sigma^2}
\right),
\end{align}{where $\sigma$ is the kernel bandwidth that controls the locality of the interpolation. We choose 
$\sigma$ proportional to the average spacing between sampled reaction-coordinate points: $\sigma = \frac{1}{d-1}\sum_{p=1}^{d-1}\Vert s_{i_{p+1}} - s_{i_{p}} \Vert$}, where $d=|\Omega|$ is the number of sampled slices. This construction can be interpreted as a kernel-based manifold interpolation of lateral slices along the reaction coordinate.


\begin{algorithm}[t!]
\caption{{BMTA Algorithm}}\label{alg:alg}
\begin{algorithmic}[1]
\State {\bf Input:} $\mathcal{C} \in \RR^{m \times d \times m}$, $\mathbf{S} \in \RR^{n \times l}$, {{$\mathbf{\Psi} \in \{0,1\}^{d \times n}$}} 
\State {\bf Parameters:}  Basis order $l$, Step size $\eta$, Max iteration $T$, target rank $r_1$ and $r_{2}$
\State {\bf Initialization:} 
{{Generate each entry of $\mathcal{\hat{G}}$, $\hat{\mathbf{X}}_1$, and $\hat{\mathbf{X}}_2$ independently from $\mathcal{N}(0,1)$}}
\State { \; {\bf Basis coefficient tensor estimation}} : Obtain the closed form solution as $\mathcal{\hat{Q}} \doteq \mbox{fold}_{(2)} \left(\mathbf{\Psi}^{\dagger}\mathbf{C}_{(2)}\right) \times_{2} \mathbf{S}^{\dagger}$.
\State {\bf \;\; Interpolative tensor estimation} Design $\hat{\mathbf{\Lambda}}_{\Omega}$
\State {\bf \;\; Low rank tensor approximation}: Using $\hat{\mathcal{Q}}\times_2 \mathbf{S}$ and $\mathcal{C} \times_2 \hat{\mathbf{\Lambda}}_{\Omega}$ obtained from the previous steps, for $t\in[T]$, update $\mathcal{\hat{G}}^{t}$ and $\mathbf{\hat{X}}^{t}_{n}$, $n=1,2$ via gradient descent method where the gradients are defined in \eqref{eq:grad_G} and \eqref{eq:grad_X} respectively.\\
    \;\;\;\; Obtain $\hat{\mathcal{ G}}=\hat{\mathcal{ G}}^{T+1}$, $\hat{\mathbf{ X}}_{n}=\hat{\mathbf{ X}}_{n}^{T+1}$ for $n=1,2$
\State {\bf Output:} $\hat{\mathcal{H}}=\hat{\mathcal{G}}\times_1 \hat{\mathbf{X}}_1 \times_2 \hat{\mathbf{X}}_2 \times_3 \hat{\mathbf{X}}_1.$
\end{algorithmic}
\end{algorithm}
\vspace*{-0.1in}

\subsubsection{3. Low-rank tensor completion} 
Finally, we exploit the quasi-structured representation in \eqref{eq:polynomial} and the interpolative tensor estimates \eqref{eq:interpolatory} to
obtain a low-rank Tucker approximation of $\mathcal{H}$. In what follows, we present the algorithm for representing low-rank tensors of $\mathcal{H}$ in terms of Tucker decomposition in \eqref{eq:Tucker}, $\mathcal{H}= \mathcal{G} \times_{1} \mathbf{X}_{1} \times_{2} \mathbf{X}_{2}\times_{3} \mathbf{X}_{1},$ where $\mathcal{G}\in\RR^{r_1 \times r_2 \times r_1}$, $r_1 \leq m$ and $r_2 \leq n$, is a low rank core tensor with $\mbox{rank}(\mathcal{G}) \leq \left[r_1,r_2,r_1\right]$, and $\mathbf{X}_{1} \in \RR^{m \times r_1}$, and $\mathbf{X}_{2} \in \RR^{n \times r_2}$ are the factor matrices to be optimized. We adopt a symmetric Tucker representation along modes 1 and 3, reflecting the symmetry of the tensor slices.
Leveraging the quasi-structured prior, $\hat{\mathcal{Q}} \times_{2} \mathbf{S}$, and the interpolative structure, $\mathcal{C} \times_2 \hat{\mathbf{\Lambda}}_{\Omega}$, obtained from the previous step, we solve the following non-convex tensor regression problem
\vspace{-0.05in}
\begin{align}\label{eq:optimization}
&\argmin_{\hat{\mathcal{G}},\hat{\mathbf{X}}_{1}, \hat{\mathbf{X}}_{2}}\; 
    \alpha\;\left\Vert \mathcal{\hat{Q}} \times_{2} \mathbf{S} - \hat{\mathcal{G}} \times_{1} \mathbf{\hat{X}}_{1} \times_{2} \mathbf{\hat{X}}_{2} \times_{3} \mathbf{\hat{X}}_{1} \right\Vert_{F}^{2}\\&\;\;\;+ (1-\alpha)\;\left\Vert {\mathcal{C}}\times_{2}\mathbf{\hat{\Lambda}}_{\Omega} -\hat{\mathcal{G}} \times_{1} \mathbf{\hat{X}}_{1} \times_{2} \mathbf{\hat{X}}_{2} \times_{3} \mathbf{\hat{X}}_{1} \right\Vert_{F}^{2},
    \nonumber
\end{align} 
where $0< \alpha <1$ is the parameter that balances the contribution of the quasi-structured prior and the manifold-guided interpolation constraint.
The proposed algorithm iteratively updates the core tensor $\mathcal{G}$ and factor matrices in an alternating fashion according to the order of :$\mathcal{G}$, $\mathbf{X}_1$, $\mathbf{X}_2$, $\mathcal{G}$, $\mathbf{X}_1$, $\dots$, $\mathbf{X}_2$ via gradient descent. Define $f(\mathcal{G},\mathbf{X}_1,\mathbf{X}_2)$ as
$f(\mathcal{G},\mathbf{X}_1,\mathbf{X}_2)=\alpha\;\left\Vert \mathcal{\hat{Q}} \times_{2} \mathbf{S} - {\mathcal{G}} \times_{1} \mathbf{{X}}_{1} \times_{2} \mathbf{{X}}_{2} \times_{3} \mathbf{{X}}_{1} \right\Vert_{F}^{2}+ (1-\alpha)\left\Vert {\mathcal{C}}\times_{2}\mathbf{\hat{\Lambda}}_{\Omega} -{\mathcal{G}} \times_{1} \mathbf{{X}}_{1} \times_{2} \mathbf{{X}}_{2} \times_{3} \mathbf{{X}}_{1} \right\Vert_{F}^{2}.$

From \eqref{eq:Kronecker}, we obtain the gradient of $f(\mathcal{G},\mathbf{X}_1,\mathbf{X}_2)$ as
\begin{align}\label{eq:grad_G}
    &\triangledown f_{\mathcal{G}}(\mathcal{G},\mathbf{X}_{1},\mathbf{X}_{2}) = -2 \alpha\left(\mathcal{\hat{Q}} \times_{2} \mathbf{S}\right) \times_{1} \mathbf{{X}}_{1}^{\trasp} \times_{2} \mathbf{{X}}_{2}^{\trasp} \times_{3} \mathbf{{X}}_{1}^{\trasp} \nonumber\\&\;\;\;\;\;- 2(1-\alpha) \left(\mathcal{{C}} \times_{2} \mathbf{\hat{\Lambda}}_{\Omega}\right) \times_{1} \mathbf{{X}}_{1}^{\trasp} \times_{2} \mathbf{{X}}_{2}^{\trasp} \times_{3} \mathbf{{X}}_{1}^{\trasp} \nonumber\\&\;\;\;\;\;+ 2 \mathcal{G} \times_{1} \mathbf{{X}}_{1}^{\trasp}\mathbf{{X}}_{1} \times_{2} \mathbf{{X}}_{2}^{\trasp}\mathbf{{X}}_{2} \times_{3} \mathbf{{X}}_{1}^{\trasp}\mathbf{{X}}_{1}.
\end{align}
We then repeat the following update rule at each iteration $t \in
\left[T\right]$ until convergence 
\begin{align}\label{eq:iterate_g}
    \mathcal{\hat{G}}^{t+1} =\mathcal{\hat{G}}^{t} - \eta \; \triangledown f_{\mathcal{\hat{G}}}\left(\mathcal{\hat{G}}^{t},\mathbf{\hat{X}}_{1}^{t},\mathbf{\hat{X}}_{2}^{t}\right),
\end{align}
where $\eta$ is a learning rate.
Next, we compute the gradient of $f(\mathcal{G},\mathbf{X}_{1},\mathbf{X}_{2})$ with respect to $\mathbf{X}_{n}$ for $n=1,2$. According to \eqref{eq:Kronecker_unfolding}, we have 
 $f(\mathcal{G},\mathbf{X}_1,\mathbf{X}_2) =\alpha\left\Vert \left(\mathcal{\hat{Q}} \times_{2} \mathbf{S}\right)_{(n)}- \mathbf{X}_{n}\mathbf{G}_{(n)}\left(\otimes_{i=1, i\neq n}^{2,1} \mathbf{X}_{i}\right)^{\trasp}\right\Vert_{F}^{2}\nonumber+ (1-\alpha)\;\left\Vert \left(\mathcal{{C}} \times_{2} \mathbf{\hat{\Lambda}}_{\Omega}\right)_{(n)}-\mathbf{X}_{n}\mathbf{G}_{(n)}\left(\otimes_{i=1, i\neq n}^{2,1} \mathbf{X}_{i}\right)^{\trasp}\right\Vert_{F}^{2}.
$ Hence, the gradient of $f(\mathcal{G},\mathbf{X}_{1},\mathbf{X}_{2})$ with respect to $\mathbf{X}_{n}$, $\triangledown f_{\mathbf{X}_{n}}(\mathcal{G},\mathbf{X}_1,\mathbf{X}_2)$, is obtained as, 
\begin{align}\label{eq:grad_X}
     &-2\alpha \left(\mathcal{\hat{Q}}\times_{2}\mathbf{S}\right)_{(n)}\mathbf{B}_{n}^{\trasp} -2(1-\alpha)\left(\mathcal{{C}}\times_{2} \mathbf{\hat{\Lambda}}_{\Omega}\right)_{(n)}\mathbf{B}_{n}^{\trasp} \nonumber\\&- 2\; \mathbf{X}_{n}\left(\mathbf{B}_{n}\mathbf{B}_{n}^{\trasp}\right), \; \mbox{where}\; \mathbf{B}_{n} = \mathbf{G}_{(n)}\left(\otimes_{i=1, i\neq n}^{2,1} \mathbf{X}_{i}\right)^{\trasp}
\end{align} Here, we do not explicitly reformulate the Kronecker product in $\mathbf{B}_{n}$, but let $\mathcal{T}=\mathcal{G} \times_1 \mathbf{X}_{1} \times_{2} \mathbf{X}_{2} \times_{3} \mathbf{X}_{1}$. Then, we have,
$\mathbf{B}_{n}=\mathbf{T}_{(n)}$, {where $\mathbf{T}_{(n)}$ is the mode-$n$ unfolding of $\mathcal{T}$ .} Finally, we update $\mathbf{X}_{n}$, for $n=1,2$ at each iteration $t \in [T]$ as follows, 
\begin{align}\label{eq:iterates}
    \mathbf{\hat{X}}_{1}^{t+1} &= \mathbf{\hat{X}}_{1}^{t} - \eta\; \triangledown f_{\mathbf{\hat{X}}_1}\left(\mathcal{\hat{G}}^{t+1},\mathbf{\hat{X}}_{1}^{t},\mathbf{\hat{X}}_{2}^{t}\right)\\
    \mathbf{\hat{X}}_{2}^{t+1} &= \mathbf{\hat{X}}_{2}^{t} - \eta\; \triangledown f_{\mathbf{\hat{X}}_2}\left(\mathcal{\hat{G}}^{t+1},\mathbf{\hat{X}}_{1}^{t+1},\mathbf{\hat{X}}_{2}^{t}\right)
\end{align}
Let $\mathcal{\hat{G}} \doteq\mathcal{\hat{G}}^{T+1}$, $\mathbf{\hat{X}}_{1} \doteq\mathbf{\hat{X}}_{1}^{T+1}$, and $\mathbf{\hat{X}}_{2} \doteq \mathbf{\hat{X}}_{2}^{T+1}$. This yields the reconstructed low-rank tensor, $\mathcal{\hat{H}}=\mathcal{\hat{G}} \times_1 \mathbf{\hat{X}}_{1} \times_2 \mathbf{\hat{X}}_{2} \times_3 \mathbf{\hat{X}}_{1}$. The proposed algorithm is summarized in Algorithm~\ref{alg:alg}.
\vspace{-0.2in}
\section{Main Theorem and Analysis}
In this section, we provide our main result and the proof sketch. The full proof can be found in the Appendix. 
\vspace{-0.15in}
\subsection{Preliminaries} We first establish that the iterates are invariant under invertible transformations of the Tucker factors. Specifically, suppose
\[
\bar{\mathbf X}_1^t=\hat{\mathbf X}_1^t\mathbf Q_1,\;\;
\bar{\mathbf X}_2^t=\hat{\mathbf X}_2^t\mathbf Q_2,\;\;
\bar{\mathcal G}^t=
(\mathbf Q_1^{-1},\mathbf Q_2^{-1},\mathbf Q_1^{-1})\cdot\hat{\mathcal G}^t,
\]
for some $\mathbf Q_k\in GL(r_k)$. Substituting these relations into \eqref{eq:iterates} shows that
\begin{eqnarray}
\bar{\mathbf X}_1^{t+1}&=\hat{\mathbf X}_1^{t+1}\mathbf Q_1, \; \;\;
\bar{\mathbf X}_2^{t+1}=\hat{\mathbf X}_2^{t+1}\mathbf Q_2,\\
\bar{\mathcal G}^{t+1}
&=
(\mathbf Q_1^{-1},\mathbf Q_2^{-1},\mathbf Q_1^{-1})\cdot\hat{\mathcal G}^{t+1}.
\end{eqnarray}
Hence, equivalent Tucker factorizations generate identical tensor iterates throughout the optimization.
\begin{align}
    \mathcal{H}^t = \left(\hat{\mathbf{X}}_{1}^t,\hat{\mathbf{X}}_{2}^t,\hat{\mathbf{X}}_{1}^t\right)\cdot\hat{\mathcal{G}}^t = \left(\bar{\mathbf{X}}_{1}^t,\bar{\mathbf{X}}_{2}^t,\bar{\mathbf{X}}_{1}^t\right)\cdot\bar{\mathcal{G}}^t,
\end{align} regardless of the representation of the tensor factors with respect to the underlying symmetry group.

\subsubsection{A distance metric} To monitor the progress of BMTA along its optimization trajectory, we require a distance metric that accounts for the inherent factor indeterminacy under invertible transformations. We therefore define the distance between factor triples $\hat{\mathbf{F}}=\left(\hat{\mathbf{X}}_1,\hat{\mathbf{X}}_2,\hat{\mathcal{G}}\right)$ and $\mathbf{F}=\left(\mathbf{X}_1,\mathbf{X}_2,\mathcal{G}\right)$ as
\begin{align}\label{eq:dist_f_fhat}
    &\mbox{dist}^{2}\left(\hat{\mathbf{F}},\mathbf{F}\right) \doteq\mbox{inf}_{\mathbf{Q}_{k}\in \mbox{GL}(r_{k})}\; \left\Vert\hat{\mathbf{X}}_{1}\mathbf{Q}_{1}-\mathbf{X}_{1}\right\Vert^{2}_{F}\nonumber \\
    &\qquad
    + \left\Vert\hat{\mathbf{X}}_{2}\mathbf{Q}_{2}-\mathbf{X}_{2}\right\Vert^{2}_{F}
    + \left\Vert\left(\mathbf{Q}_{1}^{-1},\mathbf{Q}_{2}^{-1}\right)\cdot\hat{\mathcal{G}}-\mathcal{G}\right\Vert_{F}^{2},
\end{align}
where the matrices $\{\mathbf{Q}_{k}\}_{k=1,2}$ attaining the infimum are referred to as the \emph{optimal alignment matrices}. In particular, $\hat{\mathcal{H}}$ and $\mathcal{H}$ are said to be \emph{aligned} when the optimal alignment matrices are the identity. We fix the ground-truth factors such that $\mathbf{X}_{1}$ and $\mathbf{X}_{2}$ are orthonormal matrices consisting of the left singular vectors of each mode, while the core tensor satisfies
$\mathbf{G}_{(k)}\mathbf{G}_{(k)}^{\top}=\mathbf{\Sigma}_{k}^{2},\qquad k=1,2,$
where
$
\mathbf{\Sigma}_{k}\doteq
\operatorname{diag}\!\left(
\sigma_{1}(\mathbf{H}_{(k)}),\ldots,\sigma_{r_k}(\mathbf{H}_{(k)})
\right).$
\vspace{-0.1in}
\subsection{Main result}
Denote $\sigma_{\text{min}}\left(\mathcal{H}\right) = \text{min}_{k=1,2} \;\sigma_{\text{min}}\left(\mathbf{H}_{(k)}\right)$
 as the minimum nonzero singular values of $\mathcal{H}$. In short, $\sigma_{\text{min}}$. Denote $\mathcal{E}_{\hat{\mathcal{Q}}\mathbf{S}}=\mathcal{H}-\hat{\mathcal{Q}}\times_2\mathbf{S}$ and $\mathcal{E}_{{\mathcal{C}}\hat{\mathbf{\Lambda}}_{\Omega}}=\mathcal{H}-{\mathcal{C}}\times_2\hat{\mathbf{\Lambda}}_{\Omega}.$.
We now present our main results.
\begin{theorem} \label{thm:main} Suppose that the initialization satisfies $\text{dist}\left(\hat{\mathbf{F}}_0,\mathbf{F}\right) \leq \epsilon_1 \sigma_{\min}$ for some small constant $0<\epsilon_1<1$. With probability $1-\delta$ and $d\geq O\left(\frac{n^2 l}{nl + \epsilon_2\delta(n-1)}\right)$, $0<\epsilon_2 <1$, for all $t>0$ and $\frac{1-1/\sqrt{2}}{2\text{min}_{k}\sigma_{min}^2(\mathbf{\Sigma}_k)}< \eta< \text{max}\left\{\frac{1}{2\text{min}_{k}\sigma_{max}^2(\mathbf{\Sigma}_k)},\frac{1}{1+\gamma^2(1+\tau)^{10}}\right\}$, the iterates in {\eqref{eq:grad_G} and \eqref{eq:iterates}} satisfy
\vspace{-0.05in}
{\begin{align}
&
\left\|
\left(\hat{\mathbf X}_1^{\,t},
 \hat{\mathbf X}_2^{\,t},
 \hat{\mathbf X}_1^{\,t}\right)
\cdot
\hat{\mathcal G}^{\,t}
-
\mathcal H
\right\|_F
\nonumber\\
&\le
\frac{4}{3}
\left(1+\epsilon_1\sigma_{\min}\right)^3
\zeta
\rho^{t/2}
\text{dist}\left(\hat{\mathbf{F}}_0,\mathbf F\right)
\nonumber\\
&+
O\!\left(
\frac{
\left(1+\epsilon_1\sigma_{\min}\right)^3
\zeta
}{
\sqrt{1-\rho}
}
\left[
\epsilon_1
+
\left\|
\alpha\mathbf E_{\widehat Q\mathbf S,(1)}
+
\beta\mathbf E_{\mathbf C\widehat{\Lambda}_{\Omega},(1)}
\right\|_F
\right]
\right), \nonumber
\end{align}} where $0<\rho<1$ is a constant, $\gamma=\max\left\{\sigma_1\left(\mathbf{H}_{(1)}\right),\sigma_1\left(\mathbf{H}_{(2)}\right)\right\}$, $\zeta=\max\left\{\sigma_1\left(\mathbf{H}_{(1)}\right),\sigma_1\left(\mathbf{H}_{(2)}\right),1\right\}$ and $\alpha,\beta\geq 0$, and $\alpha+\beta=1$. $\mathbf E_{\widehat Q\mathbf S,(1)}$ and $\mathbf E_{\mathbf C\widehat{\Lambda}_{\Omega},(1)}$ are the mode-1 unfolding of tensor $\mathcal{E}_{\hat{\mathcal{Q}}\mathbf{S}}$ and  $\mathcal{E}_{{\mathcal{C}}\hat{\mathbf{\Lambda}}_{\Omega
}}$. 
\end{theorem}

The Proof of Theorem ~\ref{thm:main} is in Section ~\ref{sec:mainthm}.

\vspace{-0.1in}
\subsection{Discussion}
\subsubsection{Interpreting the error} 
Theorem~\ref{thm:main} decomposes the reconstruction error into two components: a contracting optimization error and an irreducible bias due to model mismatch and sampling. The first term, proportional to $\rho^t$, decays geometrically with the optimization iterations and, by Lemma~\ref{lem:tentofac}, translates directly from the factor space to the tensor space. Consequently, when initialized within the basin of attraction, the iterates converge linearly to a neighborhood of the ground truth.
The second term depends on the quasi-structured prior approximation error $\mathcal{E}_{\hat{\mathcal{Q}}\mathbf{S}}$ and the interpolation error $\mathcal{E}_{{\mathcal{C}}\hat{\mathbf{\Lambda}}_{\Omega}}$, which quantify the mismatch between the true tensor and the assumed structural priors. Hence, the asymptotic reconstruction accuracy is determined by the quality of the basis functions and interpolation model. Theorem~\ref{thm:main} therefore characterizes the trade-off between optimization convergence and approximation accuracy.

\subsubsection{Sampling complexity $d$ and Distortion} The sampling complexity requirement on $d$ ensures that the random slice sampling operator $\mathbf{\Psi}$ satisfies a subspace embedding property with respect to the column space of the basis matrix $\mathbf{S}$, as formalized in Lemma ~\ref{lem:error}. In particular, the bound in the lemma captures the distortion induced by projecting onto the sampled slices. As $d$ increases, the embedding becomes closer to an isometry, reducing the sketching error and improving the accuracy of the basis coefficient estimate $\hat{\mathcal{Q}}$. Conversely, when $d$ is small, the deviation from isometry leads to amplification of the perturbation term $\mathcal{E}_{\hat{\mathcal{Q}}\mathbf{S}}$  which propagates through the contraction recursion in Lemma~\ref{lem:contraction}. This highlights that $d$ controls the bias–variance trade-off: larger $d$ reduces bias due to projection error but increases sampling cost. The dependence of $d$ on $n$, $l$ and  $\delta$ reflects the intrinsic dimension of the basis subspace and the concentration behavior of the sampling operator.

\subsubsection{Regularity and Closeness of $\mathbf{S}$} The known basis matrix $\mathbf{S}\in\mathbb{R}^{n\times l}$ plays a dual role in both approximation and stability. First, its approximation capability determines the magnitude of the residual tensor $\mathcal{E}_{\hat{\mathcal{Q}}\mathbf{S}}$, which directly affects the reconstruction error bound. Second, its spectral properties, including $\sigma(\mathbf{S})$, $|\mathbf{S}|_F$, and the implicit coherence in Lemma~\ref{lem:error}, govern the conditioning of the regression step \eqref{eq:Qhat} and the robustness of recovery under partial observations. Poor conditioning or high coherence in $\mathbf{S}$ amplifies the perturbation term, degrading reconstruction accuracy even when the optimization converges. In contrast, when $\mathbf{S}$ is well-conditioned and closely aligned with the true slice trajectory, the approximation error is reduced and the sampling operator preserves the relevant subspace more effectively.
Consequently, the performance of BMTA is governed by the alignment between the chosen basis representation and the intrinsic geometry of the tensor, together with the spectral gap of $\mathcal{H}$, which controls the contraction rate and stability of factor recovery.

\vspace{-0.1in}
\section{Proof and Key Lemmas for Main Result}
\label{sec:mainthm}

We provide the underlying lemmas needed to prove Theorem~\ref{thm:main} which occurs in the sequel. The complete proofs of the lemmas are deferred to the Appendix. The analysis combines a contraction argument for the gradient descent iterates in the factor space with perturbation bounds arising from sampling and model mismatch, using the alignment-invariant distance in \eqref{eq:dist_f_fhat} between the estimated factors $\hat{\mathbf{F}}_t=(\hat{\mathbf{X}}_1^t,\hat{\mathbf{X}}_2^t,\hat{\mathcal{G}}^t)$ and the ground truth $\mathbf{F}=(\mathbf{X}_1,\mathbf{X}_2,\mathcal{G})$.
Lemma~\ref{lem:contraction} establishes the contraction
\[
\mathrm{dist}^2(\hat{\mathbf{F}}_{t+1},\mathbf{F})
\le
\rho\,\mathrm{dist}^2(\hat{\mathbf{F}}_{t},\mathbf{F})
+\mathcal{E}_t,
\]
where $0<\rho<1$ and $\mathcal{E}_t$ captures the basis approximation and interpolation errors. Lemma~\ref{lem:tentofac} transfers this contraction from the factor space to the tensor space, while Lemma~\ref{lem:error} bounds the sampling-induced perturbation. Together, these results establish Theorem~\ref{thm:main}.

\begin{lemma}\label{lem:contraction} Assume that $\text{dist}\left(\hat{\mathbf{F}}_{{t}},\mathbf{F}\right) \leq \epsilon_1 \sigma_{\min}$ {for $t \geq 0$} and  for some sufficiently small $0<\epsilon_1<1$. With $\frac{1-1/\sqrt{2}}{2\text{min}_{k}\sigma_{min}^2(\mathbf{\Sigma}_k)}< \eta< \text{max}\left\{\frac{1}{2\text{min}_{k}\sigma_{max}^2(\mathbf{\Sigma}_k)},\frac{1}{1+\gamma^2(1+\tau)^{10}}\right\}$, $\gamma (1+\epsilon_1 \sigma_{min})\leq 1$, where $\gamma=\max\left\{\sigma_1\left(\mathbf{H}_{(1)}\right),\sigma_1\left(\mathbf{H}_{(2)}\right)\right\}$ and $0<\rho<1$, we have
    \begin{align}
    \text{dist}^{2}\left(\hat{\mathbf{F}}_{t+1},\mathbf{F}\right) 
    &\leq \rho\cdot\text{dist}^2\left(\hat{\mathbf{F}}_{t},\mathbf{F}\right)\nonumber \\
    &\;\;\; + O\left(\epsilon^2_1+\left\Vert\alpha\;\mathbf{E}_{\qhats,(1)} +\beta\;\mathbf{E}_{\clambda,(1)}\right\Vert_F^2\right)\nonumber
\end{align} where $\alpha+\beta=1$. 
\end{lemma}
The proof of Lemma~\ref{lem:contraction} is provided
in Appendix~\ref{prof:lem1} 

{\noindent \bf Remark. }  Lemma~\ref{lem:contraction} shows that the factor estimation error contracts linearly as long as the initialization lies within a basin of attraction. The contraction constant $\rho$ depends on the conditioning of $\mathcal{H}$ and the magnitude of the approximation error. Moreover, the convergence rate is jointly determined by the underlying tensor spectral structure $\gamma$ and the quality of the side information.

\begin{lemma}\label{lem:tentofac} Denote $\zeta=\max\left\{\sigma_1\left(\mathbf{H}_{(1)}\right),\sigma_1\left(\mathbf{H}_{(2)}\right),1\right\}$. Then, we have 
    \begin{align}
    \left\Vert\left(\hat{{\mathbf{X}}}_{1},\hat{{\mathbf{X}}}_{2},\hat{\mathbf{X}}_{1}\right)\cdot{\hat{\mathcal{G}}}-\mathcal{H}\right\Vert_{F}\leq \frac{4}{3}\left(1+\epsilon_1\sigma_{\min}\right)^3\zeta\cdot\text{dist}\left(\mathbf{F},\hat{\mathbf{F}}\right),\nonumber
\end{align} 
\end{lemma}
 The proof of Lemma~\ref{lem:tentofac} is provided
in Appendix~\ref{prof:lem2} 

{\noindent \bf Remark. }Lemma~\ref{lem:tentofac} provides a lifting argument from factor space to tensor space, showing that accurate recovery of the factors is sufficient to guarantee accurate tensor reconstruction. This step is critical since the optimization in \eqref{eq:iterate_g} and \eqref{eq:iterates} is performed over factors rather than directly over the tensor.

\begin{lemma}\label{lem:error}  {Denote the singular value decomposition (SVD) of $\mathbf{S} \in \mathbb{R}^{n \times l}$ as $\mathbf{U}_{S}\mathbf{\Sigma}_{S}\mathbf{V}_{S}^{\top}$, where $\mathbf{U}_{S}\in\mathbb{R}^{n \times l}$, $\mathbf{\Sigma}_{S}\in \mathbb{R}^{l \times l}$ and $\mathbf{V}_{S}^{\top} \in \mathbb{R}^{l\times n}$, where $n\geq l$. Assume that $\sigma_{\min}^{2}\left(\mathbf{\Psi}\mathbf{U}_{S}\right) \geq 1/\sqrt{2}$,} and  $d \geq O\left(\frac{n^2l}{nl+\delta\epsilon_2(n-1)}\right)$, where $0<\epsilon_2<1$. Then, with probability $1-\delta$, we have 
    \begin{align}
        \left\Vert\mathbf{E}_{\qhats,(1)}\right\Vert_{F}^2 \leq \left(1+2\epsilon_2\right)\left(\kappa^{-2}-1\right)\sigma_{max}^{2}\left(\mathbf{Q}_{(2)}^{*}\right)\left\Vert\mathbf{S}\right\Vert_{F}^{2},\nonumber
    \end{align} {where $\kappa \doteq
\frac{
\left\|
\mathbf U_S\mathbf U_S^\top\mathbf H_{(2)}
\right\|_F
}{
\left\|
\mathbf H_{(2)}
\right\|_F
},\;\;
\kappa\in(0,1]$ is a subspace alignment constant.} 
\end{lemma}
 The proof of Lemma~\ref{lem:error} is provided
in Appendix~\ref{prof:lem3}. 

{\noindent \bf Remark.} Lemma~\ref{lem:error} establishes a high-probability bound on the sketching error under random slice sampling, which governs the perturbation term in the contraction recursion. The bound follows from the subspace embedding concentration inequality in Lemma~\ref{lem:mat_approx2} for the randomized sampling operator $\mathbf{\Psi}$, implying that the deviation decreases as the number of observed slices increases. Moreover, the sampling complexity on $d$ ensures that the induced operator satisfies a near-isometry property over the column space of $\mathbf{S}$ \cite{dong2023fast}. Consequently, the coherence and approximation capability of $\mathbf{S}$ determine how effectively the underlying structure is preserved under partial observations.

\noindent{\bf {Proof of Theorem~\ref{thm:main}} }
Recall that Lemma~\ref{lem:contraction} establishes a contraction inequality for the distance measure,
showing that $\text{dist}^2\left(\hat{\mathbf F}_{t+1},\mathbf F\right)$
is bounded by a contraction factor $C<1$ times
$\text{dist}^2\left(\hat{\mathbf F}_{t},\mathbf F\right)$
plus an error term.
Lemma~\ref{lem:tentofac} relates the reconstruction error
$\left\Vert\left(\hat{\mathbf X}_1^t,\hat{\mathbf X}_2^t,\hat{\mathbf X}_1^t\right)\cdot \hat{\mathcal G}^t - \mathcal H\right\Vert_F$
to the distance measure $\text{dist}\left(\hat{\mathbf F}_{t},\mathbf F\right)$,
thereby controlling the optimization error in terms of the alignment distance.
Lemma~\ref{lem:error} provides an further upper bound on the sketching error
$\left\Vert\mathbf{E}_{\hat{ \mathbf{Q}}S,(1)}\right\Vert_F^2$ and sampling complexity for the number of lateral slices to obtain the corresponding error bound. Substituting the bounds from Lemmas~\ref{lem:tentofac} and \ref{lem:error} into the contraction relation of Lemma~\ref{lem:contraction} yields the recursive inequality of the form in Theorem~\ref{thm:main}. \qed
\vspace{-0.15in}
\section{Numerical Results}
In this section, we evaluate the performance of BMTA on
both synthetic and real-world datasets. All experiments on syn-
thetic data are averaged over 300 independent iterations.

\begin{table}[t]
\caption{Comparison of benchmark methods.}
\vspace*{-0.05in}
\label{tab:comparison}
\centering
\renewcommand{\arraystretch}{1.15}
\begin{tabular}{lcccc}
\toprule
Method &
Sampling &
Basis Info &
Interp. &
Low-rank \\
\midrule
Structure-only    & Lateral  & \checkmark & & \checkmark \\
Interp-only & Lateral  & & \checkmark & \checkmark \\
Tucker             & Adapted lateral  & & & \checkmark \\
TensorCUR          & Lateral + Horizontal  & & & \checkmark \\
\textbf{BMTA}      & Lateral  & \checkmark & \checkmark & \checkmark \\
\bottomrule
\end{tabular}
\vspace*{-0.15in}
\end{table}
\subsubsection{Benchmark algorithms}
Table~\ref{tab:comparison} summarizes the methods considered in our experiments. We compare BMTA with two ablation variants, {\em Structure-only} and {\em Interpolation-only}, which exploit only the quasi-structure model in \eqref{eq:polynomial} and the interpolation model in \eqref{eq:interpolatory}, respectively. We also compare against Tucker decomposition \cite{malik2018low,kolda2009tensor} and TensorCUR \cite{mahoney2006tensor}. Unless otherwise specified, all methods use the same sampling budget $d$.
Since Tucker decomposition requires a fully observed tensor, we first construct an initial estimate
\begin{equation}
\tilde{\mathcal H}(:,j,:)
=
\frac{1}{d}
\sum_{i\in\Omega}
\mathcal H(:,i,:),
\qquad
j\notin\Omega,
\end{equation}
where $\Omega$ denotes the sampled lateral slice indices. Tucker decomposition is then obtained by solving
$
\min_{\mathcal G,X_1,X_2}
\left\|
\tilde{\mathcal H}
-
\mathcal G
\times_1X_1
\times_2X_2
\times_3X_1
\right\|_F^2
$
via gradient descent.
TensorCUR reconstructs the tensor from sampled lateral and horizontal slices using the CUR decomposition $\hat{\mathcal H}=\mathcal C\mathcal U\mathcal R,$
where $\mathcal U$ is computed from the intersection tensor using the Moore--Penrose pseudoinverse. Unlike BMTA, TensorCUR exploits only sampled observations and does not incorporate structural or interpolation priors.
\vspace{-0.1in}
\subsection{Synthetic datasets}
\subsubsection{Data generation} 
The ground-truth tensor $\mathcal{H} \in \mathbb{R}^{m \times n \times m}$ is constructed as a superposition of two components. {For our synthetic datasets in this section, we use a Legendre polynomial basis for $\mathbf{S}$.}

Let $\mathbf{S} \in \mathbb{R}^{n \times l}$ denote a polynomial basis matrix of degree $l$, and let $\mathcal{Q} \in \mathbb{R}^{m \times l \times m}$ denote a coefficient tensor. The entries of $\mathcal{Q}$ are defined as
$
\mathcal{Q}(i,k,j) = \frac{1}{(k+1)^{\gamma}} \sum_{m=1}^{K} a_m 
\exp\!\left(-\frac{(x_i - c_{x,m})^2}{2\sigma_{x,m}^2} - \frac{(y_j - c_{y,m})^2}{2\sigma_{y,m}^2}\right),\nonumber
$
where $a_m$, $c_{x,m}$, $c_{y,m}$, 
$\sigma_{x,m}$, and $\sigma_{y,m}$ are randomly generated. And we obtain polynomial structure via the mode-2 tensor-matrix product
$\mathcal{H}_{\text{poly}} = \mathcal{Q} \times_2 \mathbf{S},$
i.e., $\mathcal{H}_{\text{poly}}(i,t,j) = \sum_{k=1}^{l} \mathcal{Q}(i,k,j)\, \mathbf{S}(t,k).$
The interpolative component is defined as a time-varying term:
\begin{equation}
\mathcal{H}_{\text{inter}}(i,t,j) = -\exp\!\left(-\frac{(x_i - c_x(t))^2 + (y_j - c_y(t))^2}{2\sigma^2}\right), \label{eq:interp_num}
\end{equation}
where $(c_x(t),c_y(t))$ denotes the time-varying center location of 
the local perturbation, which evolves smoothly along the trajectory. Finally, both components are normalized as
$\widetilde{\mathcal{H}}_{\text{poly}} = \frac{\mathcal{H}_{\text{poly}}}{\|\mathcal{H}_{\text{poly}}\|_F}, 
\quad
\widetilde{\mathcal{H}}_{\text{local}} = \frac{\mathcal{H}_{\text{inter}}}{\|\mathcal{H}_{\text{inter}}\|_F}.$
The final tensor is constructed as
\begin{equation}\label{eq:true_synthetic}
\mathcal{H} = (1-\rho)\,\widetilde{\mathcal{H}}_{\text{poly}} 
+ \rho\,\widetilde{\mathcal{H}}_{\text{local}} + \mathcal{E},
\end{equation}
where $\rho \in [0,1]$ controls the relative contribution of the two components, and $\mathcal{E}$ is a noise tensor with Gaussian  normal distribution.

\subsubsection{Varying $d$} 
\begin{figure}[t]
    \centering
    \includegraphics[width=0.85\linewidth]{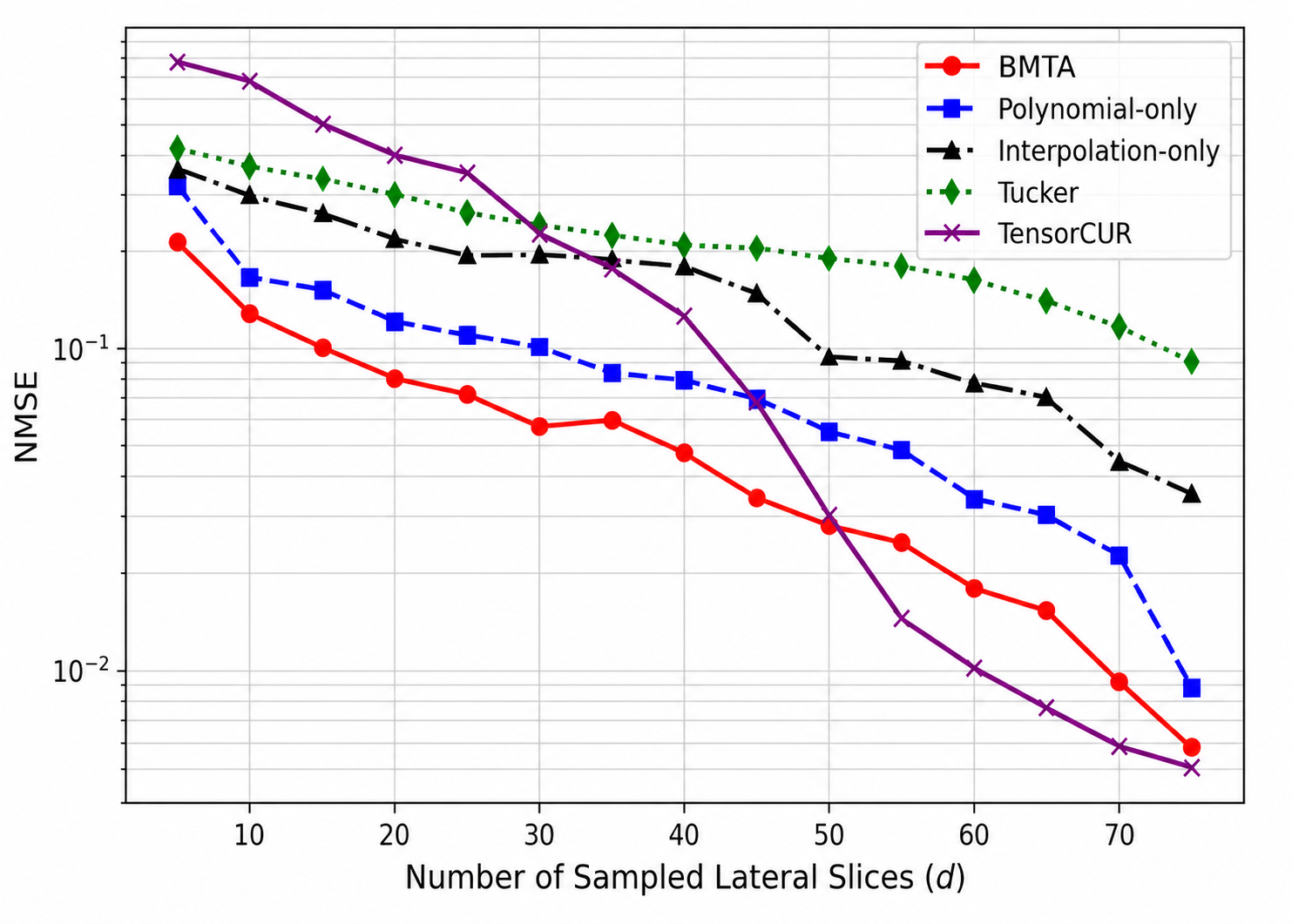}
    \vspace*{-0.15in}
    \caption{NMSE comparison varying the number of sampled lateral slices $d$}
    \label{fig:d_nmse}
\end{figure}

Figure~\ref{fig:d_nmse} shows the reconstruction NMSE as a function of the number of sampled lateral slices 
$d$, where $\alpha=0.5$, and $r_1=r_2=5$. In the low-to-moderate sampling regime, BMTA consistently achieves the lowest error, reflecting its ability to incorporate structured side information to regularize the reconstruction under limited observations. In contrast, Polynomial-only and Interpolation-only rely on fixed priors and cannot fully capture the underlying variability, while Tucker decomposition suffers from unreliable subspace estimation when the number of sampled slices is small.{ The results closely follow the trends predicted by Theorem~\ref{thm:main}. As the number of sampled lateral slices $d$ increases, the sampling-induced perturbation decreases (more accurate subspace estimation in step \eqref{eq:Qhat}), resulting in improved reconstruction accuracy (lower $\left\|
\alpha\mathbf E_{\widehat Q\mathbf S,(1)}
+
\beta\mathbf E_{\mathbf C\widehat{\Lambda}_{\Omega},(1)}
\right\|_F$) across all methods. } TensorCUR eventually outperforms others in the high-sampling regime, where data-driven recovery becomes reliable. {However, it must be observed that TensorCUR has access to the true temporally evolving nature of the signal via the horizontal slices, whereas BMTA conjectures the interpolative model which is mismatched with respect to the data generation in \eqref{eq:interp_num}.}

\begin{figure}[t]
    \centering
    \includegraphics[width=0.85\linewidth]{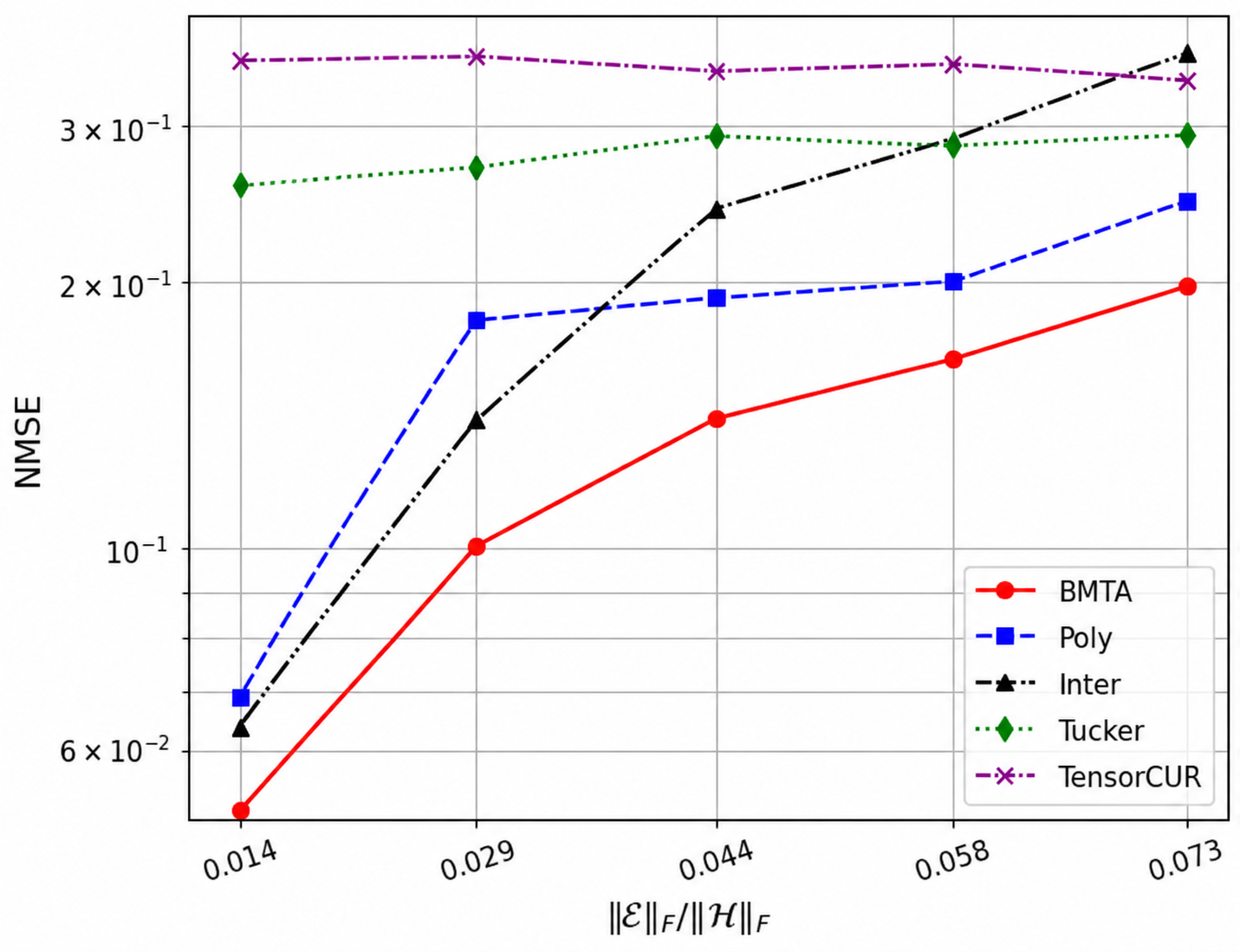
    }
        \vspace*{-0.15in}
    \caption{NMSE versus normalized noise level}
    \label{fig:noise}
        \vspace*{-0.2in}
\end{figure}

\begin{figure*}[t]
    \centering
\includegraphics[width=15.5cm]{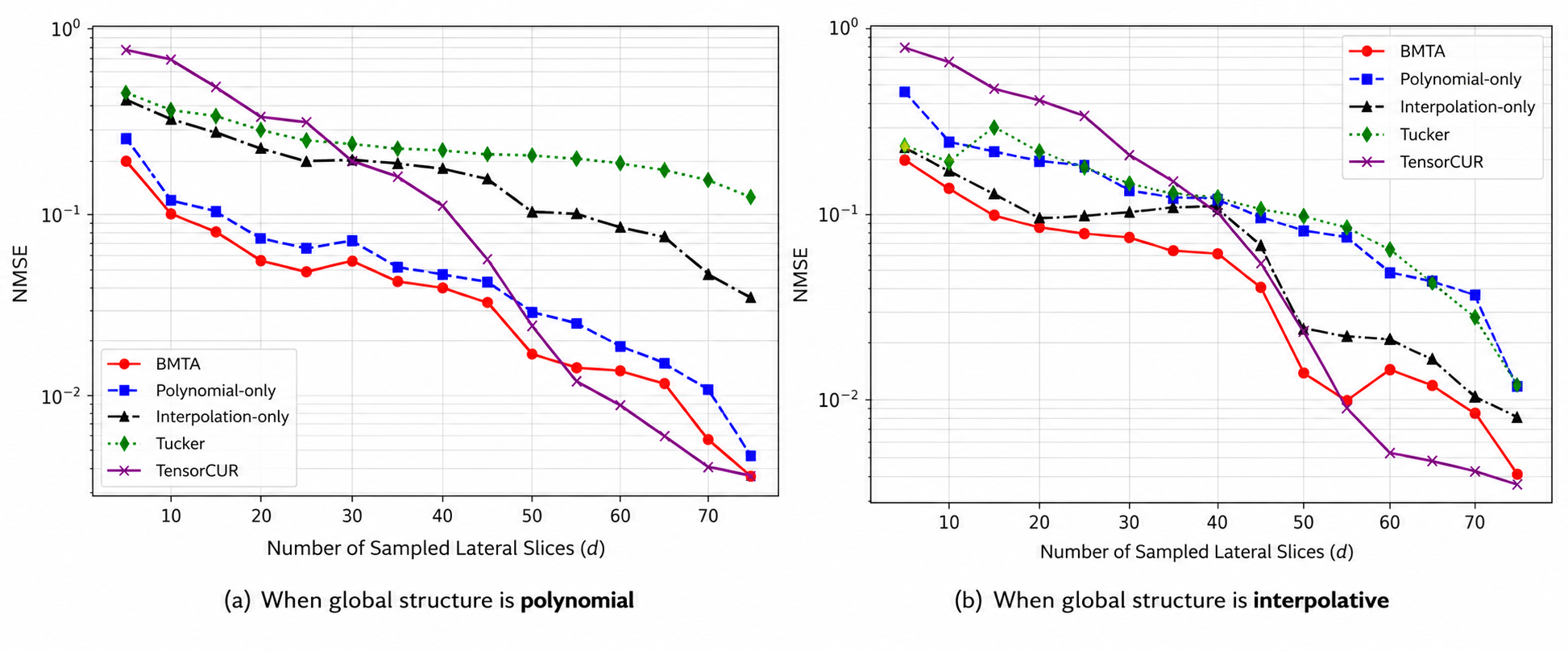}
    \vspace{-0.2in}
    \caption{Sensitivity of BMTA to the underlying tensor structure. NMSE comparison when the tensor is dominated by (a) quasi-polynomial structure and (b) interpolative structure.}
    \label{fig:global_structure}
    \vspace{-0.25in}
\end{figure*}

\subsubsection{Sensitivity to noise} We next investigate the robustness of BMTA to noisy observations in Fig.~\ref{fig:noise}. We generate the true tensor varying the normalized noise level $(\Vert\mathcal{E}\Vert_F/\Vert\mathcal{H}\Vert_F)$. As the noise level increases, the NMSE of all algorithms increases due to the larger perturbation from the underlying structured tensor model. {This observation is consistent with Theorem~\ref{thm:main}, where the reconstruction error depends on the mismatch between the true tensor and the imposed quasi-structured prior and interpolative structures.}
We observe that BMTA consistently achieves the lowest NMSE across all noise levels. This is because BMTA does not rely only on the sampled tensor slices, but instead combines two complementary structural priors: the quasi-polynomial trajectory and manifold-guided interpolation. The additional structural information regularizes the reconstruction problem and improves robustness against noisy observations. In contrast, while Tucker and TensorCUR rely primarily on sampled tensor observations and therefore cannot exploit the underlying tensor evolution.

\subsubsection{Sensitivity to the underlying tensor structure} We next evaluate how the performance of BMTA depends on the dominant structure of the tensor in Fig.~\ref{fig:global_structure}. Recall that the synthetic tensor is generated as a combination of a quasi-polynomial component and an interpolative component in \eqref{eq:true_synthetic}, where $\rho$ controls their relative contributions.
Fig.~\ref{fig:global_structure}(a) shows the quasi-polynomial-dominant regime. As expected, Polynomial-only achieves competitive performance since the assumed model closely matches the true tensor generation process. Nevertheless, BMTA further improves the reconstruction accuracy by incorporating the interpolation prior, which captures local deviations beyond the global polynomial basis.
Similarly, Fig.~\ref{fig:global_structure}(b) shows the interpolation-dominant regime. Although Interpolation-only benefits from the strong local correlation among lateral slices, BMTA consistently achieves lower NMSE by additionally exploiting the global structural prior. This observation aligns with Theorem~\ref{thm:main}, which shows that the reconstruction error depends jointly on both structural approximation terms, explaining the robustness of BMTA through the complementary use of the structured basis and manifold-guided interpolation priors.

\vspace{-0.1in}
\subsection{Application: Radio map reconstruction}

\textbf{Data generation.}
We evaluate BMTA on a synthetic spatio-temporal radio map reconstruction task motivated by wireless network planning and base-station placement \cite{yi2026stfd,luo2025denoising}. We extend the radio map generation model in \cite{9693274,romero2022radio} to the tensor setting by generating a sequence of time-varying radio maps.
The radio environment is defined on an $80\times80$ spatial grid over $[0,L]^2$ with $L=2000$~m, yielding a tensor $\mathcal{H}\in\mathbb{R}^{80\times80\times80}$. The received power is generated using the log-distance path-loss model
$
P_{\mathrm{PL}}(s)
=
10\log_{10}
\left(
\sum_{k=1}^{3}
10^{P_k(s)/10}
\right),$
where
$
P_k(s)
=
P_{0,k}
-
10n_k
\log_{10}
\left(
\frac{\max(\|s-s_k\|,d_0)}{d_0}
\right),
$
with $d_0=10$~m, $P_{0,k}\in\{-30,-35,-38\}$~dBm, and $n_k\in\{2.3,2.6,2.1\}$.

The global temporal evolution is modeled by a DCT basis,
$
\mathcal{H}_{\mathrm{DCT}}
=
\mathcal{Q}\times_2\mathbf{S},$
while local variations are modeled by a smoothly varying shadowing field
$
\mathcal{H}_{\mathrm{inter}}(x,t,y)
=
\sum_{i=1}^{10}
a_i(t)
\cos
\left(
2\pi
\frac{k_{x,i}x+k_{y,i}y}{L}
+
\phi_i(t)
\right),$
where the amplitudes and phases vary smoothly with time. Both components are normalized before combination, and the final tensor is generated as
$
\mathcal{H}
=
(1-\rho)\tilde{\mathcal{H}}_{\mathrm{DCT}}
+
\rho\tilde{\mathcal{H}}_{\mathrm{inter}}
+
\mathcal{E},
$
where $\rho=0.5$ unless otherwise specified and $\mathcal{E}$ is additive Gaussian noise. The number of observed lateral slices varies from $d=10$ to $70$.

\subsubsection{Performance analysis} 
\begin{figure}[t]
    \centering
    \includegraphics[width=0.75\linewidth]{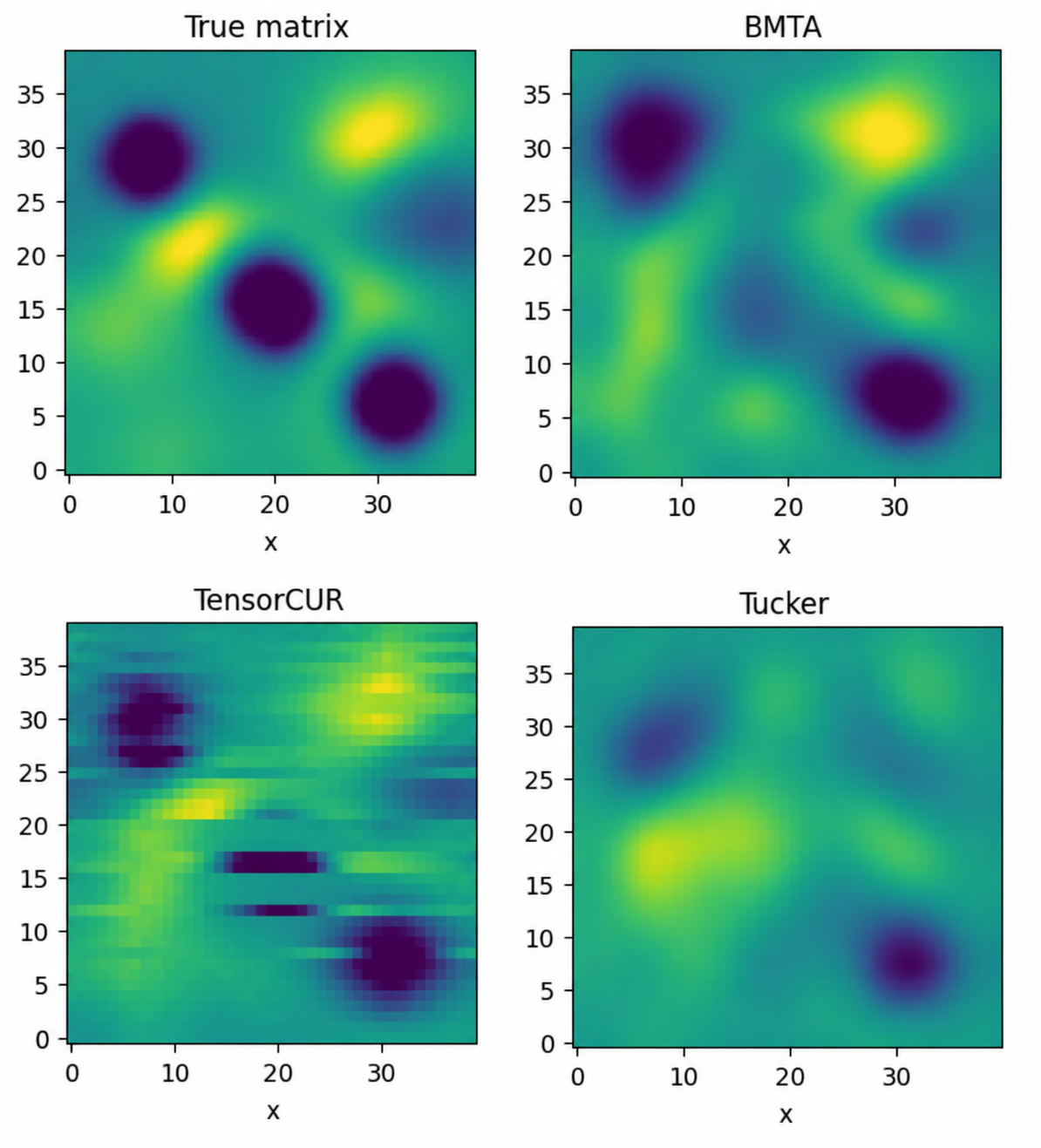}
    \vspace{-0.15in}
    \caption{Comparison on radio map dataset at $t = 65$: BMTA provides the superior reconstruction of spatio-temporal features.}
    \label{fig:radiomap}
    \vspace*{-0.2in}
\end{figure}

Table~\ref{tb:radio_map} shows the reconstruction performance on the synthetic radio map dataset. BMTA achieves the lowest reconstruction error in most sampling regimes, especially when the number of observed slices is limited (modest-sampling regime). {When the number of slices is large ($d=60,70$), TensorCUR benefits from access to the true interpolative nature of the signal via horizontal slices.} Our numerical results demonstrates that the proposed combination of DCT-based side information and interpolative modeling in BMTA effectively compensates for missing tensor observations when the sampled observation is limited. As shown in Fig.~\ref{fig:radiomap}, BMTA better preserves the spatial patterns of the true radio map, while TensorCUR and Tucker exhibit artifacts or overly smoothed reconstructions.
\vspace{-0.1in}
\begin{table}[t]
\centering
\caption{Performance comparison on synthetic radio map dataset}
\begin{tabular}{c|ccc}
\toprule
$d$ & {\bf BMTA} & TensorCUR \cite{mahoney2006tensor} & Tucker \cite{kolda2009tensor} \\
\midrule
10 & 0.4727 & 0.7586 & 0.5674 \\
20 & 0.2965 & 0.4705 & 0.4605 \\
30 & 0.2456 & 0.3126 & 0.3973 \\
40 & 0.1511 & 0.1524 & 0.3245 \\
50 & 0.1012 & 0.1073 & 0.2547 \\
60 & 0.0981 & 0.0671 & 0.1961 \\
70 & 0.0622 & 0.0371 & 0.1028 \\
\bottomrule
\end{tabular}\label{tb:radio_map}
\vspace*{-0.15in}
\end{table}

\begin{table}[t!]
\centering
\caption{Performance comparison for multiple real-world reactions {{when $d=5$}}}
    \begin{center}
    \begin{tabular}{c|c|c|c|c|c}
    \toprule
    System & BMTA & Poly.& Interp. & Tucker & TensorCUR  \\  \hline
    CF3CH3  & {\bf 0.0254}  & 0.048   & 0.042 & 0.087  & 0.0361\\ \hline
    TSoxo & {\bf 0.0318} & 0.0424 & 0.0322 & 0.065 & 0.0338\\ \hline
    Sn2Ar1 & {\bf 0.0053} & 0.0068 & 0.0063 & 0.013 & 0.0095\\\hline
    Ir & {\bf 0.0271} & 0.0601& 0.0333 & 0.085 & 0.0523\\ \hline
    \bottomrule
    \end{tabular}
    \label{tb:quantum_chem}
    \end{center}
    \vspace{-0.2in}
\end{table} 
\vspace{-0.1in}
\subsection{Application: Quantum Chemistry}
{
Table~\ref{tb:quantum_chem} compares the reconstruction performance on multiple real-world chemical reaction datasets \cite{bac2022matrix,bac2025incorporating}. {Based on our prior knowledge on the reaction systems, we use polynomial basis for $\mathbf{S}$. } BMTA achieves the lowest reconstruction error across all systems by integrating polynomial trajectory and interpolative relationships among tensor slices. We use $l=3$ for system CF3CH3 and Sn2Ar1, and $l=5$ for TSoxo and Ir depending on the complexity of the system.  Although Polynomial-only and Interpolation-only can capture either global or local structures, they show the limited reconstruction accuracy. 
BMTA exhibits superior performance versus TensorCUR, which further demonstrates the benefit of exploiting domain-specific side information under limited slice observations.}

\vspace{-0.25in}
\section{Conclusions}

In this work, we proposed Basis and Manifold prior Tensor Approximation (BMTA), a tensor approximation framework for scenarios where only a limited number of tensor slices are available. BMTA exploits two complementary sources of side information: a basis representation that captures the global tensor evolution and a manifold-guided interpolation model that characterizes local relationships among lateral slices. By incorporating these priors into a low-rank Tucker optimization framework, BMTA enables accurate reconstruction from highly constrained slice observations. We established theoretical guarantees through a reconstruction error bound that characterizes the convergence behavior, sampling complexity, and approximation errors arising from structural model mismatch. Numerical experiments on synthetic and real-world datasets demonstrated that BMTA consistently improves reconstruction accuracy over existing tensor approximation methods, particularly in limited-sampling regimes. These results highlight the effectiveness of integrating complementary structural priors for recovering high-dimensional tensor data from limited observations.
\appendices
\vspace{-0.2in}
\section{Technical preliminaries}
\vspace*{-0.1in}
In what follows, we provide several useful lemmas that are used to prove our main theorem and additional lemmas.
\begin{lemma}\label{lem:mat_approx2}
    \cite{drineas2006fast} Let $\mathbf{A}$ be a $n \times m$ matrix, and $\mathbf{B}$ be a $n \times p$ matrix. Let $s$ denote the number of samples. Draw $s$ i.i.d row indices uniformly at random without replacement, and form the matrix $\left(\mathbf{\Psi A}\right)^{\top}\mathbf{\Psi B}$, where $\mathbf{\Psi}$ is the random {$n \times n$} row sampling and rescaling operator {which has only $s$ non-zero components by construction.}  {the matrix $\left(\mathbf{\Psi A}\right)^{\top}\mathbf{\Psi B}$, approximates $\mathbf{A}^{\top}\mathbf{B}$}.  With probability at least $1-\delta$, the approximation satisfies,\vspace{-0.05in}
    \begin{align}
    &\left\Vert\mathbf{A}^{\top}\mathbf{B}-\left(\mathbf{\Psi A}\right)^{\top}\mathbf{\Psi B}\right\Vert_{F}^{2}\nonumber\\&\;\;\;\; \leq \frac{n}{(n-1)\delta}\left(\frac{n}{s}-1\right)\sum_{k=1}^{n}\vert\mathbf{A}(k,:)\vert^{2}\vert\mathbf{B}(k,:)\vert^{2}.\nonumber
    \end{align}
\end{lemma}
\begin{lemma}\label{lem:mat_approx}
    {Denote the singular value decomposition (SVD) of $\mathbf{S} \in \mathbb{R}^{n \times l}$ as $\mathbf{U}_{S}\mathbf{\Sigma}_{S}\mathbf{V}_{S}^{\top}$, where $\mathbf{U}_{S}\in\mathbb{R}^{n \times l}$, $\mathbf{\Sigma}_{S}\in \mathbb{R}^{l \times l}$ and $\mathbf{V}_{S}^{\top} \in \mathbb{R}^{l\times n}$, where $n\geq l$.} {Denote the mode-2 unfolding of
a tensor $\mathcal{H}$ as $\mathbf{H}_{(2)}$. Then, we define
$\mathbf{H}_{(2)}^\perp\doteq\left(\mathbf{I}-\mathbf{U}_S\mathbf{U}_S^\top\right)\mathbf{H}_{(2)}$.}  Construct a random slice-sampling matrix $\mathbf{\Psi}$ that samples each row of $\mathbf{S}$ uniformly at random without replacement. Then, provided $d\geq O\left(\frac{n^2 l}{nl + \epsilon_2\delta(n-1)}\right)$, where $\epsilon_2>0$, then, $\left\Vert\mathbf{U}_{S}^{\top}\mathbf{\Psi}^{\top}\mathbf{\Psi}\mathbf{H}_{(2)}^{\perp}\right\Vert_{F}^{2} \leq \epsilon_2\left\Vert\mathbf{H}_{(2)}^{\perp}\right\Vert^{2}$ holds with probability $\delta$.
\end{lemma}


{\noindent\em{Proof of Lemma~\ref{lem:mat_approx}.}}  
We can apply Lemma~\ref{lem:mat_approx2} to show that  Lemma~\ref{lem:mat_approx} holds,
\begin{align} \mathbb{E}\left[\left\Vert\mathbf{U}_{S}^{\top}\mathbf{\Psi}^{\top}\mathbf{\Psi}\mathbf{H}_{(2)}^{\perp}\right\Vert_{F}^{2}\right]   &\stackrel{{(a)}}{=} \mathbb{E}\left[\left\Vert\mathbf{U}_{S}^{\top}\mathbf{H}_{(2)}^{\perp}-\mathbf{U}_{S}^{\top}\mathbf{\Psi}^{\top}\mathbf{\Psi}\mathbf{H}_{(2)}^{\perp}\right\Vert_{F}^{2}\right]\nonumber\\
    &\stackrel{(b)}{\leq}\frac{nl}{\delta(n-1)}\left(\frac{n}{d}-1\right)\left\Vert\mathbf{H}_{(2)}^{\perp}\right\Vert_{F}^{2},\label{eq:lem2_1}
\end{align}    \vspace*{-0.015in} \noindent {where $(a)$ holds $\mathbf{U}_{S}^{\top}\mathbf{H}_{(2)}^{\perp} = \mathbf{U}_{S}^{\top}\left(\mathbf{I}-\mathbf{U}_S\mathbf{U}_S^\top\right)\mathbf{H}_{(2)}=\left(\mathbf{U}_{S}^{\top}-\mathbf{U}_{S}^{\top}\right)\mathbf{H}_{(2)}=0$
as $\mathbf{H}_{(2)}^\perp$ represents the component of $\mathbf{H}_{(2)}$ orthogonal to the column space of $\mathbf{S}$.}  $(b)$ is due to Lemma~\ref{lem:mat_approx2} and the fact that $\mathbf{U}_S$ is unitary.
For non-negative random variable $\mathbf{X}$ and scalar $t > 0$, the Markov inequality states that $\mathbf{X} \geq t$ as $\mathbb{P}\left[\mathbf{X} \geq t\right] \leq \mathbb{E}[\mathbf{X}]/t$.  We apply the Markov inequality to bound the probability that the slice sampling matrix violates Lemma~\ref{lem:mat_approx}:
\begin{align}
    &\mathbb{P}_{\mathbf{S}\sim\mathcal{D}}\left[\left\Vert\mathbf{U}_{S}^{\top}\mathbf{\Psi}^{\top}\mathbf{\Psi}\mathbf{H}_{(2)}^{\perp}\right\Vert_{F}^{2}\geq \epsilon_2\left\Vert\mathbf{H}_{(2)}^{\perp}\right\Vert_{F}^{2}\right]\nonumber\\&\;\;\leq \frac{\mathbb{E}\left[\left\Vert\mathbf{U}_{S}^{\top}\mathbf{\Psi}^{\top}\mathbf{\Psi}\mathbf{H}_{(2)}^{\perp}\right\Vert_{F}^{2}\right]}{\epsilon_2\left\Vert\mathbf{H}_{(2)}^{\perp}\right\Vert_{F}^{2}} \stackrel{(a)}{\leq} \frac{nl}{\delta\epsilon_2(n-1)}\left(\frac{n}{d}-1\right), 
\end{align} \noindent where $(a)$ is from \eqref{eq:lem2_1}. Therefore, when $d \geq O\left(\frac{n^2l}{nl+\delta\epsilon_2(n-1)}\right)$,  with probability $1-\delta$, we have $\left\Vert\mathbf{U}_{S}^{\top}\mathbf{\Psi}^{\top}\mathbf{\Psi}\mathbf{H}_{(2)}^{\perp}\right\Vert_{F}^{2} \leq \epsilon_2\left\Vert\mathbf{H}_{(2)}^{\perp}\right\Vert_{F}^{2}$.\qed
\vspace{-0.1in}
\section{Several perturbation bounds}
We will repeatedly use the following relationships, which follow directly from standard matrix relations applied to the matricized forms \cite{kolda2009tensor}. {Let $\mathbf{A}_k \in \mathbb{R}^{r_k \times r_k}$ for $k=1,2,3$, and $\mathcal{S} \in \mathbb{R}^{r_1 \times r_2 \times r_3}$ is a tensor.} 
For a tensor $\mathcal{X}$ of multilinear rank at most
$r = (r_1,r_2,r_3)$,
its spectral norm and Frobenius norm satisfy
\cite{tong2022scaling} $\|\mathcal{X}\|_F
\le
\sqrt{\frac{r_1 r_2 r_3}{r}}\,\|\mathcal{X}\|$ where $r = \max_{k=1,2,3} r_k.$  Thus, we have,
\begin{align}
(\mathbf{X}_1,\mathbf{X}_2,\mathbf{X}_3)\big((\mathbf{A}_1,\mathbf{A}_2,\mathbf{A}_3)\cdot \mathcal{S}\big)
&=
\left(\mathbf{X}_1\mathbf{A}_1,\, \mathbf{X}_2\mathbf{A}_2,\, \mathbf{X}_3\mathbf{A}_3\right)\cdot \mathcal{S}\nonumber\\
\|\left(\mathbf{A}_1,\mathbf{A}_2,\mathbf{A}_3\right)\cdot \mathcal{S}\|_F
&\leq
\|\mathbf{A}_1\|_2\,\|\mathbf{A}_2\|_2\,\|\mathbf{A}_3\|_2\,\|\mathcal{S}\|_F,\nonumber
\end{align}

We introduce several perturbation bounds. First, we introduce the following notation that will be used repeatedly:\vspace{
-0.05in
}
\begin{align}
    &\Delta_{\hat{\mathbf{X}}_1}\doteq \hat{\mathbf{X}}_{1}-\mathbf{X}_{1}, \;\; \Delta_{\hat{\mathbf{X}}_2}\doteq \hat{\mathbf{X}}_{2}-\mathbf{X}_{2},\;\;\Delta_{\hat{\mathcal{G}}}\doteq \hat{\mathcal{G}}-\mathcal{G} \\&\xcheck\doteq \label{eq:not1}
    \left(\hat{\mathbf{X}}_2\otimes\hat{\mathbf{X}}_1\right)\hat{\mathbf{G}}^{\trasp}_{(1)}, \;\;\;\; \check{\mathbf{X}}_{2}\doteq 
    \left(\hat{\mathbf{X}}_2\otimes\hat{\mathbf{X}}_1\right)\hat{\mathbf{G}}^{\trasp}_{(2)}\\
    &\bar{\mathbf{X}}_{1}\doteq
    \left({\mathbf{X}}_2\otimes{\mathbf{X}}_1\right){\mathbf{G}}^{\trasp}_{(1)}, \;\;\;\; \bar{\mathbf{X}}_{2}\doteq 
    \left({\mathbf{X}}_2\otimes{\mathbf{X}}_1\right){\mathbf{G}}^{\trasp}_{(2)}.
\end{align} {Recall $\mathbf{G}_{(k)}$ denotes the mode-$k$ matricization of $\mathcal{G}$.} {Let $\vee$ denote the maximum operation.} Now we are ready to state the lemma on perturbation bounds.
\vspace{-0.1cm}
\begin{lemma}\label{lem:spec_bound} Suppose $\mathbf{F}=\left(\mathbf{X}_1,\mathbf{X}_2,\mathcal{G}\right)$ and $\hat{\mathbf{F}}=\left(\hat{\mathbf{X}}_1,\hat{\mathbf{X}}_2,\hat{\mathcal{G}}\right)$ are aligned and satisfy $\text{dist}\left(\hat{\mathbf{F}},\mathbf{F}\right)\leq \epsilon_1\sigma_{\min}$ for some $\epsilon_1 < 1$. Denote $\Delta_{\hat{\mathbf{G}},(k)}$ as mode-$k$ unfolding of $\Delta_{\hat{\mathcal{G}}}$. Then the following bounds hold regarding the
spectral norm:
\vspace{-0.1in}
\begin{align}
\left\Vert\Delta_{\hat{\mathbf{X}}_{1}}\right\Vert \vee \left\Vert\Delta_{\hat{\mathbf{X}}_{2}}\right\Vert \vee\left\Vert{\Delta_{\hat{\mathbf{G}},(1)}}\right\Vert &\leq \epsilon_1\sigma_{min}\label{eq:b1}\\
\left\Vert{\mathbf{\Sigma}_{1}^{-1}\Delta_{\hat{\mathbf{G}},(1)}}\right\Vert&\leq \epsilon_1\label{eq:b2}\\
    \left\Vert\left(\xcheck-\xbar\right)\mathbf{\Sigma}_{1}^{-1}\right\Vert\leq 3\epsilon_1+4\epsilon_1^2&\sigma_{min}+2\epsilon_1^3\sigma^2_{min}\label{eq:b3}
\end{align}
\end{lemma}
{\noindent\em{Proof of Lemma~\ref{lem:spec_bound}.}} 
{Since $\mathbf{F}$ and $\hat{\mathbf{F}}$ are aligned, the optimal
alignment matrices in the definition of
$\operatorname{dist}(\hat{\mathbf{F}},\mathbf{F})$
are the identity matrices. Therefore,
\vspace{-0.05in}
\begin{align}
\operatorname{dist}\!\left(\hat{\mathbf{F}},\mathbf{F}\right)
=\sqrt{
\left\|\Delta_{\hat{\mathbf{X}}_1}\right\|_{F}^{2}
+
\left\|\Delta_{\hat{\mathbf{X}}_2}\right\|_{F}^{2}
+
\left\|\Delta_{\hat{\mathcal{G}}}\right\|_{F}^{2}}.\label{eq:dist_def}
\end{align} 
\vspace{-0.05in}
By the assumption
$
\text{dist}\!\left(\hat{\mathbf{F}},\mathbf{F}\right)
\leq
\epsilon_1\sigma_{\min},
$
each nonnegative term is individually bounded by the total distance.
Since the spectral norm is bounded above by the Frobenius norm, we
obtain
$
\left\|\Delta_{\hat{\mathbf{X}}_1}\right\|
\vee
\left\|\Delta_{\hat{\mathbf{X}}_2}\right\| \vee \left\| \Delta_{\hat{\mathbf{G}},(1)}\right\|
\leq
\epsilon_1\sigma_{\min},
$ which proves \eqref{eq:b1}.
Next, for each $k=1,2$, submultiplicativity of the spectral norm gives $\left\|
\mathbf{\Sigma}_k^{-1}
\Delta_{\hat{\mathbf{G}},{(k)}}
\right\|_2
\leq
\left\|\mathbf{\Sigma}_k^{-1}\right\|_2
\left\|\Delta_{\hat{\mathbf{G}},{(k)}}\right\|_F
\leq
\frac{
\left\|\Delta_{\hat{\mathbf{G}},{(k)}}\right\|_{F}
}{
\sigma_{\min}(\mathbf{\Sigma}_k)
}.$
By the definition
$
\sigma_{\min}(\mathcal{H})
=
\min_{j=1,2}\sigma_{\min}(\mathbf{\Sigma}_j),
$
we have $\sigma_{\min}(\mathbf{\Sigma}_k)
\geq
\sigma_{\min}(\mathcal{H}).$
It follows that
$\left\|
\mathbf{\Sigma}_k^{-1}
\Delta_{\hat{\mathbf{G}}_{(k)}}
\right\|_2\leq
\frac{
\epsilon_1\sigma_{\min}(\mathcal{H})
}{
\sigma_{\min}(\mathbf{\Sigma}_k)
}\leq
\epsilon_1,$
which proves \eqref{eq:b2}. 
Equation \eqref{eq:b3} is shown next.
Observe that
$
\check{\mathbf X}_1-\bar{\mathbf X}_1
=
\Big[
(\hat{\mathbf X}_2\otimes\hat{\mathbf X}_1)
-
(\mathbf X_2\otimes\mathbf X_1)
\Big]
\mathbf G_{(1)}^\top
+
(\hat{\mathbf X}_2\otimes\hat{\mathbf X}_1)
\Delta_{\hat{\mathbf G},(1)}^\top .
$
Multiplying by $\mathbf\Sigma_1^{-1}$ and using the triangle inequality yields,
$
\|
(\check{\mathbf X}_1-\bar{\mathbf X}_1)\mathbf\Sigma_1^{-1}
\|_2
\le
T_1+T_2,
$
where $T_1=
\Big\|
\big[
(\hat{\mathbf X}_2\otimes\hat{\mathbf X}_1)
-
(\mathbf X_2\otimes\mathbf X_1)
\big]
\mathbf G_{(1)}^\top\mathbf\Sigma_1^{-1}
\Big\|,$ and $T_2=
\|
(\hat{\mathbf X}_2\otimes\hat{\mathbf X}_1)
\Delta_{\hat{\mathbf G},(1)}^\top
\mathbf\Sigma_1^{-1}
\|.$ We have,
{\begin{align}
    T_1 &= \Big\|
\big[
(\hat{\mathbf X}_2\otimes\hat{\mathbf X}_1)
-
(\mathbf X_2\otimes\mathbf X_1)
\big]
\mathbf G_{(1)}^\top\mathbf\Sigma_1^{-1}
\Big\|_2\nonumber\\
&\stackrel{(a)}{\leq} \left\Vert\Delta_{\hat{\mathbf{X}}_2}\otimes\mathbf{X}_1+\mathbf{X}_2 \otimes \Delta_{\hat{\mathbf{X}}_1} + \Delta_{\hat{\mathbf{X}}_2} \otimes \Delta_{\hat{\mathbf{X}}_1} \right\Vert_2 \left\Vert\mathbf G_{(1)}^\top\mathbf\Sigma_1^{-1}\right\Vert_2\nonumber\\
&\stackrel{(b)}{\leq} \left(\left\Vert\Delta_{\hat{\mathbf{X}}_2}\right\Vert_2+\left\Vert\Delta_{\hat{\mathbf{X}}_1}\right\Vert_2+\left\Vert\Delta_{\hat{\mathbf{X}}_1}\right\Vert_2\left\Vert\Delta_{\hat{\mathbf{X}}_2}\right\Vert_2\right)\nonumber\\
&\stackrel{(c)}{\leq} 2\epsilon_1\sigma_{\min}
+
\epsilon_1^2\sigma_{\min}^2,
\end{align}}

\noindent where $(a)$ follows from submultiplicativity of norms, 
$(b)$  follows from the facts that $\mathbf G_{(1)}\mathbf G_{(1)}^\top=\mathbf\Sigma_1^2$, and
$\|\mathbf G_{(1)}^\top\mathbf\Sigma_1^{-1}\|_2=1$,
$\hat{\mathbf X}_k \doteq \mathbf X_k+\Delta_{\hat{\mathbf X}_k}$,
$\|\mathbf X_k\|=1$,
and
$\|\mathbf A\otimes\mathbf B\|_2
=\|\mathbf A\|_2\,\|\mathbf B\|_2$; and $(c)$ follows from \eqref{eq:b1}.
Similarly,
$T_2 \leq
\|\hat{\mathbf X}_2\|_2
\|\hat{\mathbf X}_1\|_2
\|
\Delta_{\hat{\mathbf G},(1)}^\top
\mathbf\Sigma_1^{-1}
\|_2
\le
\epsilon_1
(1+\epsilon_1\sigma_{\min})^2,$
where we use
$\|\hat{\mathbf X}_k\|_2
\le
1+\epsilon_1\sigma_{\min}$.
Combining the above bounds completes the proof of \eqref{eq:b3}. \qed}

Next, we first recall a result on the existence of optimal
alignment matrices for Tucker factorizations which will be used in this paper.
{\begin{lemma}\cite{dong2023fast}
\label{lem:scaled_inf}
Let $\mathbf{F}
=\left(\mathbf{X}_1,\mathbf{X}_2,\mathbf{X}_3,\mathcal{G}\right)$ and
$\hat{\mathbf{F}}
=
\left(\hat{\mathbf{X}}_1,\hat{\mathbf{X}}_2,\hat{\mathbf{X}}_3,\hat{\mathcal{G}}\right),$
and define the scaled distance
$
\text{dist}_s^2\left(\hat{\mathbf F},\mathbf F\right)
\doteq
\inf_{\mathbf Q_k\in\mathrm{GL}(r_k)}
\sum_{k=1}^{3}
\left\|
\left(\hat{\mathbf X}_k\mathbf Q_k-\mathbf X_k\right)
\mathbf\Sigma_k
\right\|_F^2
+\left\|
\left(\mathbf Q_1^{-1},
 \mathbf Q_2^{-1},
 \mathbf Q_3^{-1}\right)
\cdot \hat{\mathcal G}-\mathcal G
\right\|_F^2,$
where $\mathbf\Sigma_k
=
\text{diag}
\!\left(
\sigma_1(\mathcal H_{(k)}),
\dots,
\sigma_{r_k}(\mathcal H_{(k)})
\right),\;\;
k=1,2,3.$
If $\text{dist}_s\left(\hat{\mathbf F},\mathbf F\right)
<
\sigma_{\min}\left(\mathcal H\right),$ then the infimum is attained by some $\mathbf Q_k\in\mbox{GL}(r_k)$ for
$k=1,2,3.$
\end{lemma}
We now specialize this result to the symmetric Tucker decomposition considered in this paper.
\begin{lemma}\label{lem:inf}
    Define the distance measure between $\hat{\mathbf{F}}=\left(\hat{\mathbf{X}}_1, \hat{\mathbf{X}}_2, \hat{\mathcal{G}}\right)$ and ${\mathbf{F}}=\left({\mathbf{X}}_1, {\mathbf{X}}_2, {\mathcal{G}}\right)$ as 
\begin{align}
    &\mbox{dist}^{2}\left(\hat{\mathbf{F}},\mathbf{F}\right) \doteq \mbox{inf}_{\mathbf{Q},\mathbf{R}\in \mbox{GL}(r_{k})}\; \left\Vert\left(\hat{\mathbf{X}}_{1}\mathbf{Q}-\mathbf{X}_{1}\right)\right\Vert_{F}^{2}\nonumber \\&+ \left\Vert\left(\hat{\mathbf{X}}_{2}\mathbf{R}-\mathbf{X}_{2}\right)\right\Vert_{F}^{2} + \left\Vert\left(\mathbf{Q}^{-1},\mathbf{R}^{-1},\mathbf{Q}^{-1}\right)\cdot\hat{\mathcal{G}}-\mathcal{G}\right\Vert_{F}^{2}\nonumber.
\end{align}  
Then, suppose $\mbox{dist}\left(\hat{\mathbf{F}},\mathbf{F}\right)<\sigma_{min}$, then there exist some $\mathbf{Q} \in GL(r_1)$ and $\mathbf{R} \in GL(r_2)$, which attain the infimum of the above distance.
\end{lemma}}



\vspace{-0.3in}
\section{Proof of Lemma~\ref{lem:contraction}}\label{prof:lem1} 
\vspace*{-0.01in}
{Recall that $\mathcal{E}_{\hat{\mathcal{Q}}\mathbf{S}}=\mathcal{H}-\hat{\mathcal{Q}}\times_2\mathbf{S}$ and $\mathcal{E}_{{\mathcal{C}}\hat{\mathbf{\Lambda}}_{\Omega}}=\mathcal{H}-{\mathcal{C}}\times_2\hat{\mathbf{\Lambda}}_{\Omega}.$}
{From our assumptions in the lemma,} we have for {$t\geq 0$},
$
\text{dist}(\hat{\mathbf F}_t,\mathbf F)
\le
\epsilon_1\sigma_{\min}$ for
$0<\epsilon_1<1.$
Let $\mathbf Q_t\in\mathrm{GL}(r_1)$ and
$\mathbf R_t\in\mathrm{GL}(r_2)$ denote the optimal alignment
matrices between $\hat{\mathbf F}_t$ and $\mathbf F$, whose
existence is guaranteed by Lemma~\ref{lem:inf}. For notational
simplicity, we absorb these alignment matrices into the current
iterate by defining
$\hat{\mathbf X}_1
\doteq
\hat{\mathbf X}_{t,1}\mathbf Q_t, 
\hat{\mathbf X}_2
\doteq
\hat{\mathbf X}_{t,2}\mathbf R_t,$
and
$\hat{\mathcal G}
\doteq
(\mathbf Q_t^{-1},\mathbf R_t^{-1},\mathbf Q_t^{-1})
\cdot
\hat{\mathcal G}_t.$
Accordingly, we write
$
\hat{\mathbf F}
\doteq
(\hat{\mathbf X}_1,\hat{\mathbf X}_2,\hat{\mathcal G}).$
{
To bound
$\text{dist}(\hat{\mathbf F}_{t+1},\mathbf F)$,
we use the same alignment matrices
$\mathbf Q_t$ and $\mathbf R_t$.
Although
$\mathbf Q_t$
and
$\mathbf R_t$
are optimal for
$\hat{\mathbf F}_t$,
they are not necessarily optimal for
$\hat{\mathbf F}_{t+1}$.
Nevertheless, they provide a feasible candidate in the
minimization defining
$\text{dist}(\hat{\mathbf F}_{t+1},\mathbf F)$. Therefore, {from the definition of the distance measure in Equation~(\ref{eq:dist_f_fhat}), we have}
\begin{align}\label{eq:11}
\text{dist}^{2}
(\hat{\mathbf F}_{t+1},\mathbf F)
&\le
\left\|
\hat{\mathbf X}_{t+1,1}\mathbf Q_t-\mathbf X_1
\right\|_F^2
+
\left\|
\hat{\mathbf X}_{t+1,2}\mathbf R_t-\mathbf X_2
\right\|_F^2
\nonumber\\
&+
\left\|
(\mathbf Q_t^{-1},\mathbf R_t^{-1},\mathbf Q_t^{-1})
\cdot
\hat{\mathcal G}_{t+1}
-
\mathcal G
\right\|_F^2.
\end{align}}
Our goal is to bound the distance in \eqref{eq:11} from the ground truth to the next iterate, i.e. $\text{dist}\left(\hat{\mathbf{F}}_{t+1},\mathbf{F}\right)$. We first bound $\left\|
\hat{\mathbf X}_{t+1,1}\mathbf Q_t-\mathbf X_1
\right\|_F^2$. Using the gradient descent update rule in \eqref{eq:grad_X} together with the expression
for $\nabla_{\mathbf {\hat{X}}_1}f$, we have
$
\hat{\mathbf X}_{t+1,1}\mathbf Q_t-\mathbf X_1
=
\hat{\mathbf X}_1-\mathbf X_1
-\eta\nabla_{\mathbf{\hat{X}}_1}f
=
\Delta_{\hat{\mathbf X}_1}
-
2\eta
\left[
\alpha
\left(
\hat{\mathcal H}
-
(\mathcal H+\mathcal E_{\qhats})
\right)_{(1)}
+
\beta
\left(
\hat{\mathcal H}
-
(\mathcal H+\mathcal E_{\clambda})
\right)_{(1)}
\right]
\check{\mathbf X}_1,$
where
$
\hat{\mathcal H}
\doteq
(\hat{\mathbf X}_1,\hat{\mathbf X}_2,\hat{\mathbf X}_1)
\cdot\hat{\mathcal G}$, $\Delta_{\hat{\mathbf X}_1}
\doteq
\hat{\mathbf X}_1-\mathbf X_1,$
and
$
\check{\mathbf X}_1
\doteq
(\hat{\mathbf X}_2\otimes\hat{\mathbf X}_1)
\hat{\mathbf G}_{(1)}^\top.$
Since $\alpha+\beta=1$, the bracketed term 
above can be written as
$(\hat{\mathcal H}-\mathcal H)_{(1)}
-
\alpha\mathcal E_{\qhats,(1)}
-
\beta\mathcal E_{\clambda,(1)}.$
Moreover, using
$
\hat{\mathcal H}_{(1)}
=
\hat{\mathbf X}_1\check{\mathbf X}_1^\top,$ and $
\mathcal H_{(1)}=\mathbf X_1\bar{\mathbf X}_1^\top,$
where
$\bar{\mathbf X}_1\doteq
(\mathbf X_2\otimes\mathbf X_1)
\mathbf G_{(1)}^\top,$
we obtain
\vspace{-0.05in}
\begin{align}
(\hat{\mathcal H}-\mathcal H)_{(1)}
&=
\Delta_{\hat{\mathbf X}_1}
\check{\mathbf X}_1^\top
+
\mathbf X_1
\left(
\check{\mathbf X}_1-\bar{\mathbf X}_1
\right)^\top
\label{eq:residual-decomposition}
\end{align}
\vspace{-0.05in}
{Recall that $\mathbf{E}_{\qhats,(k)}$ and $\mathbf{E}_{\clambda,(k)}$ are mode-$k$ unfolding of $\mathcal E_{\qhats}$ and $\mathcal E_{\clambda}$, respectively.}
Substituting {\eqref{eq:residual-decomposition} into $(\hat{\mathcal H}-\mathcal H)_{(1)}$} yields,
\begin{align}
&\hat{\mathbf X}_{t+1,1}\mathbf Q_t-\mathbf X_1\nonumber\\&=\Delta_{\hat{\mathbf X}_1}(\mathbf{I}-2\eta\check{\mathbf X}_1^\top\check{\mathbf X}_1) -2\eta\mathbf{X}_{1}\left(\xcheck-\xbar\right)^{\trasp}\xcheck\nonumber\\&\;\;\;\;\;\;\;\;\;+2\alpha\eta\mathbf{E}_{\qhats,(1)}\xcheck +2\beta\eta\mathbf{E}_{\clambda,(1)}\xcheck,\label{eq:diff}
\end{align}
Let $\odot  \doteq \Delta_{\hat{\mathbf X}_1}(\mathbf{I}-2\eta\check{\mathbf X}_1^\top\check{\mathbf X}_1) -2\eta\mathbf{X}_{1}\left(\xcheck-\xbar\right)^{\trasp}\xcheck$. {We next compute the squared two-norm of each side of (\ref{eq:diff}), and exploit the triangle inequality and submultiplicativity of norms,}
\begin{align}
&\left\Vert\hat{\mathbf{X}}_{t+1,1}\mathbf{Q}_{t}-\mathbf{X}_{1}\right\Vert^{2}_2\nonumber\\
&\leq \left\Vert\odot\right\Vert^{2}_2 +4\eta\left\Vert\odot\right\Vert_2\left\Vert\alpha\mathbf{E}_{\qhats,(1)}+\beta\mathbf{E}_{\clambda,(1)}\right\Vert_2\left\Vert\xcheck\right\Vert_2\nonumber\\&\;\;\;\;\;+4\eta^{2}\left\Vert\alpha\mathbf{E}_{\qhats,(1)} +\beta\mathbf{E}_{\clambda,(1)}\right\Vert^{2}_2\left\Vert\xcheck\right\Vert^2_2.\label{eq:6}
\end{align}
Using the same tools as above, we bound $ \left\Vert\odot\right\Vert^{2}_2$ as,

\vspace{-0.05in}
\begin{align}
    \left\Vert\odot\right\Vert^{2}_2
    &\leq 2\left\Vert\mathbf{I}-2\eta\check{\mathbf X}_1^\top\check{\mathbf X}_1\right\Vert^2\left\Vert\Delta_{\hat{\mathbf{X}}_{1}}\right\Vert^{2}+8\eta^{2}\underbrace{\left\Vert\mathbf{X}_{1}\left(\xcheck-\xbar\right)^{\trasp}\xcheck\right\Vert^{2}}_{\odot_1}.\nonumber
\end{align} 
We next bound $ \sqrt{\odot_{1}}$ using submultiplicativity, 
\begin{align}
    \sqrt{\odot_{1}} &= \left\Vert\mathbf{X}_{1}\left(\xcheck-\xbar\right)^{\trasp}\xcheck\right\Vert_{2}
    \stackrel{(a)}{=}\left\Vert\left(\xcheck-\xbar\right)^{\trasp}\mathbf{\Sigma}_1^{-1}\mathbf{\Sigma}_1\xcheck\right\Vert_{2}\nonumber\\
    &\leq \left\Vert\left(\xcheck-\xbar\right)\mathbf{\Sigma}_1^{-1}\right\Vert_{2}\left\Vert\mathbf{\Sigma}_1\xcheck\right\Vert_{2},\label{eq:5}
\end{align} where $(a)$ is due to the orthogonality of $\mathbf X_1$.  The term
$\left\Vert\mathbf{\Sigma}_1\xcheck\right\Vert_{2}$ is bounded as
{{Recall that $\mathbf{G}_{(k)}$ denotes the mode-$k$ matricization of $\mathcal{G}$.}}
{\begin{align}
&\left\Vert\mathbf{\Sigma}_1\xcheck\right\Vert_2\stackrel{(a)}{\leq} \left\Vert\mathbf{\Sigma}_1\right\Vert_2\left\Vert\hat{\mathbf{X}}_{2}\right\Vert_2\left\Vert\hat{\mathbf{X}}_{1}\right\Vert_2\left\Vert\mathbf{\Sigma}_1\mathbf{\Sigma}_1^{-1}\hat{\mathbf{G}}_{(1)}\right\Vert_2\nonumber\\
    &=\left\Vert\mathbf{\Sigma}_1\right\Vert_2\left\Vert\hat{\mathbf{X}}_{2}\right\Vert_2\left\Vert\hat{\mathbf{X}}_{1}\right\Vert_2\left\Vert\mathbf{\Sigma}_1\mathbf{\Sigma}_1^{-1}\left({\mathbf{G}}_{(1)}-{\Delta_{\hat{\mathcal{G}},(1)}}\right)\right\Vert_2\nonumber\\
    &\stackrel{(b)}{\leq}\left\Vert\mathbf{\Sigma}_1\right\Vert(1+\epsilon_1\sigma_{min})^{2}\left(\left\Vert\mathbf{\Sigma}_1\mathbf{\Sigma}_1^{-1}{\mathbf{G}_{(1)}}\right\Vert+\left\Vert\mathbf{\Sigma}_1\mathbf{\Sigma}_1^{-1}{\Delta_{\hat{\mathcal{G}},(1)}}\right\Vert\right)\nonumber\\
    &\stackrel{(c)}{\leq} (1+\epsilon_1)(1+\epsilon_1\sigma_{min})^2\left\Vert\mathbf{\Sigma}_1\right\Vert_2^2\label{eq:xcheck},
\end{align}}  \noindent {where $(a)$ is due to norm inequality and \eqref{eq:not1} and $(b)$ and $(c)$ are due to \eqref{eq:b1} and \eqref{eq:b2} respectively.}
Plugging \eqref{eq:xcheck} into \eqref{eq:5}, we have
\begin{align}
    \sqrt{\odot_{1}}
    &\stackrel{(a)}{\leq}(1+\epsilon_1)(1+\epsilon_1\sigma_{min})^2\left(3\epsilon_1+4\epsilon_1^2\sigma_{min}+2\epsilon_1^3\sigma^2_{min}\right)\left\Vert\mathbf{\Sigma}_{1}\right\Vert^2\nonumber\\
&\stackrel{(b)}{\leq}
6\epsilon_1(1+\epsilon_1\sigma_{min})^4\left\Vert\mathbf{\Sigma}_1\right\Vert^2
\end{align} \noindent { where $(a)$ is due to \eqref{eq:b3} and $(b)$ follows from $0<\epsilon_1\leq 1$, and the fact that if $c \doteq\epsilon_1\sigma_{\min}$, then
$3+4c+2c^2\leq 3(1+c)^2$. }
To simplify the notation, define
\begin{align}
q
&\doteq
\left\|
\alpha\mathbf E_{\qhats,(1)}
+
\beta\mathbf E_{\clambda,(1)}
\right\|,\\
K_1
&\doteq
6(1+\epsilon_1\sigma_{min})^4\left\Vert\mathbf{\Sigma}_1\right\Vert^2,\\
L_1
&\doteq(1+\epsilon_1)
\left(1+\epsilon_1\sigma_{\min}\right)^2
\|\mathbf\Sigma_1\|.
\end{align}
From the previous result,
$
\sqrt{\odot_1}\leq \epsilon_1K_1,$. Assuming $0<\eta\leq \frac{1}{2\text{max}_{k}\sigma^2_{max}(\mathbf{\Sigma}_k)}$ for $k=1,2$, we have $\left\Vert\mathbf{I}-2\eta\check{\mathbf X}_1^\top\check{\mathbf X}_1\right\Vert=1-2\eta\sigma_{min}^2(\mathbf{\Sigma}_k)$,
which implies
\begin{align}
\|\odot\|^2
&\leq
2\left(1-2\eta\sigma_{min}^2(\mathbf{\Sigma}_1)\right)^2
\|\Delta_{\hat{\mathbf X}_1}\|_F^2
+
8\eta^2\odot_1
\nonumber\\
&\stackrel{(a)}{\leq}
2\left(1-2\eta\sigma_{min}^2(\mathbf{\Sigma}_1)\right)^2
\operatorname{dist}^2
\left(\hat{\mathbf F}_t,\mathbf F\right)
+
8\eta^2\epsilon_1^2K_1^2,\nonumber
\end{align} \noindent where $(a)$ is due to \eqref{eq:dist_def}.
Recalling the bound of
$\|\hat{\mathbf X}_{t+1,1}\mathbf Q_t-\mathbf X_1\|^2$
in \eqref{eq:6}, we have
\begin{align}
&
\left\|
\hat{\mathbf X}_{t+1,1}\mathbf Q_t-\mathbf X_1
\right\|_F^2\leq
\|\odot\|^2
+
4\eta q
\|\odot\|
\|\check{\mathbf X}_1\|
+
4\eta^2q^2
\|\check{\mathbf X}_1\|^2
\nonumber\\
&\leq
2\left(1-2\eta\sigma_{min}^2(\mathbf{\Sigma}_1)\right)^2
\operatorname{dist}^2
\left(\hat{\mathbf F}_t,\mathbf F\right)
+
8\eta^2\epsilon_1^2K_1^2
\nonumber\\
&+
4\eta qL_1
\sqrt{
\left(1-2\eta\sigma_{min}^2(\mathbf{\Sigma}_1)\right)^2
\operatorname{dist}^2
\left(\hat{\mathbf F}_t,\mathbf F\right)
+
4\eta^2\epsilon_1^2K_1^2
}
\nonumber\\&\;\;\;\;\;+
4\eta^2q^2L_1^2.\nonumber
\end{align}
Now define
$
A\doteq2\left(1-2\eta\sigma_{min}^2(\mathbf{\Sigma}_1)\right)
\operatorname{dist}
\left(\hat{\mathbf F}_t,\mathbf F\right)$, $B\doteq
2\sqrt{2}\eta\epsilon_1K_1,$ and
$C\doteq 2\eta qL_1.$
Then,
\begin{align}
&
\left\|
\hat{\mathbf X}_{t+1,1}\mathbf Q_t-\mathbf X_1
\right\|_F^2\leq
A^2+B^2+2C\sqrt{A^2+B^2}+C^2
\nonumber\\
&=
\left(\sqrt{A^2+B^2}+C\right)^2
\stackrel{(a)}{\leq}
(A+B+C)^2
\nonumber\\
&\stackrel{(b)}{\leq}
4\left(1-2\eta\sigma_{min}^2(\mathbf{\Sigma}_1)\right)^2
\operatorname{dist}^2
\left(\hat{\mathbf F}_t,\mathbf F\right)+
4\eta^2
(\sqrt{2}\epsilon_1K_1+qL_1)^2,\nonumber
\end{align} 
\noindent
where $(a)$ follows from $\sqrt{x^2+y^2}\leq x+y$ for $x,y\ge0$, and
$(b)$ follows from Young's inequality
$(x+y)^2\leq 2x^2+2y^2$. For $\frac{1-1/\sqrt{2}}{2\text{min}_{k}\sigma_{min}^2(\mathbf{\Sigma}_k)}< \eta< \frac{1}{2\text{max}_{k}\sigma_{max}^2(\mathbf{\Sigma}_k)}$, $2\left(1-2\eta\sigma_{min}^2(\mathbf{\Sigma}_1)\right)^2 < 1$,
we finally obtain
\vspace{-0.05in}
\begin{align}
&
\left\|
\hat{\mathbf X}_{t+1,1}\mathbf Q_t-\mathbf X_1
\right\|_F^2
\leq
\rho_0
\operatorname{dist}^2
\left(\hat{\mathbf F}_t,\mathbf F\right)
+
8\eta^2
(\epsilon_1K_1+qL_1)^2, 
\label{eq:xq_x1}
\end{align}\noindent  where $\rho_0 <1$. {Using similar tools we determine the following bound:}
$\left\Vert\hat{\mathbf{X}}_{t+1,2}\mathbf{R}_{t}-\mathbf{X}_{2}\right\Vert_{F}^{2}\leq \rho_1\;\mbox{dist}^2\left(\hat{\mathbf{F}}_{t},\mathbf{F}\right) + 8\eta^2
(\epsilon_1K_2+qL_2)^2$, where $K_2 = 6(1+\epsilon_1\sigma_{min})^4\left\Vert\mathbf{\Sigma}_2\right\Vert^2$ and $L_2=(1+\epsilon_1)\left(1+\epsilon_1\sigma_{\min}\right)^2
\|\mathbf\Sigma_2\|$ for $\rho_1<1$.
We next bound $\left\Vert\left(\mathbf{Q}_{t}^{-1},\mathbf{R}_{t}^{-1},\mathbf{Q}_{t}^{-1}\right)\cdot\hat{\mathcal{G}}_{t+1}-\mathcal{G}\right\Vert_{F}^{2}$. Again, we use the update rule and $\hat{\mathcal{G}} = \Delta_{\mathcal{G}} + \mathcal{G}$ to obtain $\left(\mathbf{Q}_{t}^{-1},\mathbf{R}_{t}^{-1},\mathbf{Q}_{t}^{-1}\right)\cdot\hat{\mathcal{G}}_{t+1}-\mathcal{G} =(1-\eta)\Delta_{\hat{\mathcal{G}}}-\eta\left({\mathcal{G}} \times_{1} \mathbf{\hat{X}}_{1} \times_{2} \mathbf{\hat{X}}_{2} \times_{3} \mathbf{\hat{X}}_{1} -\mathbf{W}\right)\cdot \hat{\mathbf{{X}}}_{1}^{\top} \cdot \hat{\mathbf{{X}}}_{2}^{\top} \cdot \hat{\mathbf{{X}}}_{1}^{\top}$
\noindent
where $\mathbf{W}=\alpha\left(\mathcal{\hat{Q}} \times_{2} \mathbf{S}\right)-(1-\alpha)\left({\mathcal{C}}\times_{2}\mathbf{\hat{\Lambda}}_{\Omega}\right).$
Taking squared norm on each side, we have 
\vspace{-0.05in}
\begin{align} &\left\Vert(\mathbf{Q}_{t}^{-1},\mathbf{R}_{t}^{-1},\mathbf{Q}_{t}^{-1})\cdot\hat{\mathcal{G}}_{t+1}-\mathcal{G}\right\Vert_{F}^{2}\stackrel{(a)}{\leq} 2\left\Vert(1-\eta)\Delta_{\hat{\mathcal{G}}}\right\Vert_{F}^{2}\nonumber\\ 
    & + 2\eta^2\underbrace{\left\Vert\left({\mathcal{G}} \cdot \mathbf{\hat{X}}_{1} \cdot \mathbf{\hat{X}}_{2} \cdot \mathbf{\hat{X}}_{1} -\mathbf{W}\right)\cdot \mathbf{\hat{X}}_{1}^{\top} \cdot\mathbf{\hat{X}}_{2}^{\top} \cdot\mathbf{\hat{X}}_{1}^{\top}\right\Vert_{F}^{2}}_{\circledast_1}\label{eq:hatg_g},
    \vspace{-0.05in}
\end{align} 

\noindent
where $(a)$ is due to the norm triangle inequality.

\subsubsection{Bound of $\circledast_{1}$} We use ${\mathcal{G}} \cdot \mathbf{\hat{X}}_{1} \cdot \mathbf{\hat{X}}_{2} \cdot \mathbf{\hat{X}}_{1} -\mathbf{W} = {\mathcal{G}} \cdot \mathbf{\hat{X}}_{1} \cdot \mathbf{\hat{X}}_{2} \cdot \mathbf{\hat{X}}_{1} - \left(\mathcal{H} - \alpha\mathcal{E}_{\qhats}-(1-\alpha)\mathcal{E}_{\clambda}\right)
    $ 
to rewrite $\circledast_{1}$; we define $\mathbf{T} \doteq \alpha\mathcal{E}_{\qhats}-(1-\alpha)\mathcal{E}_{\clambda}$, thus,
\begin{align}
    \circledast_{1} &= \left\Vert\left({\mathcal{G}} \cdot \mathbf{\hat{X}}_{1} \cdot \mathbf{\hat{X}}_{2} \cdot \mathbf{\hat{X}}_{1} - \left(\mathcal{H} - \mathbf{T}\right)\right)\cdot \mathbf{\hat{X}}_{1}^{\top} \cdot\mathbf{\hat{X}}_{2}^{\top} \cdot\mathbf{\hat{X}}_{1}^{\top}\right\Vert_{F}^2\nonumber\\
    &\stackrel{(a)}{\leq} \left\Vert\mathbf{\hat{X}}_{1}^{\top} \cdot\mathbf{\hat{X}}_{2}^{\top} \cdot\mathbf{\hat{X}}_{1}^{\top}\right\Vert_{2}^2\left(\left\Vert{\mathcal{G}} \cdot \mathbf{\hat{X}}_{1} \cdot \mathbf{\hat{X}}_{2} \cdot \mathbf{\hat{X}}_{1} - \mathcal{H}\right\Vert_{F}^2 + \left\Vert\mathbf{T}\right\Vert_{F}^2\right)\nonumber\\
    &\stackrel{(b)}{\leq}  (1+\epsilon_1\sigma_{min})^6   \! \! \left( \!\underbrace{\left\Vert{\mathcal{G}} \cdot \mathbf{\hat{X}}_{1} \cdot \mathbf{\hat{X}}_{2} \cdot \mathbf{\hat{X}}_{1} - \mathcal{H}\right\Vert_{F}^2}_{\circledast_{1,1}}  \! \! + \left\Vert\mathbf{T}\right\Vert_{F}^2 \!\right) \label{eq:g_1},
\end{align} \noindent {where $(a)$ follows from the submultiplicativity of Frobenius and spectral norms and
$\|(\mathbf{A}-\mathbf{B})\cdot\mathbf C\|_F
\leq
\left(\|\mathbf{A}\|_F+
\|\mathbf{B}\|_F\right)
\|\mathbf C\|_2$. Inequality
$(b)$ follows from the triangle inequality  and \eqref{eq:b1} for $k=1,2$  yielding
$\|\hat{\mathbf X}_k\|
\leq
\|\mathbf X_k\|+\|\Delta_{\hat{\mathbf X}_k}\|
\leq
1+\epsilon_1\sigma_{\min}$
}

Let $\gamma=
\max\left\{
\|\mathbf{\Sigma}_1\|,
\|\mathbf{\Sigma}_2\|
\right\}$ and $\tau=\epsilon_1\sigma_{min}$. Now, we bound $\circledast_{1,1}$ as follows. 
\begin{align}
    &{\mathcal{G}} \cdot \mathbf{\hat{X}}_{1} \cdot \mathbf{\hat{X}}_{2} \cdot \mathbf{\hat{X}}_{1} - \mathcal{H}= {\mathcal{G}} \cdot \mathbf{\hat{X}}_{1} \cdot \mathbf{\hat{X}}_{2} \cdot \mathbf{\hat{X}}_{1} - {\mathcal{G}} \cdot \mathbf{{X}}_{1} \cdot \mathbf{{X}}_{2} \cdot \mathbf{{X}}_{1}\nonumber\\
    &= {\mathcal{G}} \cdot \left(\Delta_{\hat{\mathbf{X}}_1}, \hat{\mathbf{X}}_2, \hat{\mathbf{X}}_1\right) + {\mathcal{G}} \cdot \left(\hat{\mathbf{X}}_1, \Delta_{\hat{\mathbf{X}}_2}, \hat{\mathbf{X}}_1\right) + {\mathcal{G}} \cdot \left(\hat{\mathbf{X}}_1, \hat{\mathbf{X}}_2, \Delta_{\hat{\mathbf{X}}_1}\right)\nonumber
\end{align} Using this equality, we obtain,
\begin{align}
    &\left\Vert{\mathcal{G}} \cdot \mathbf{\hat{X}}_{1} \cdot \mathbf{\hat{X}}_{2} \cdot \mathbf{\hat{X}}_{1} - \mathcal{H}\right\Vert_{F}
    \nonumber\\&\stackrel{(a)}{\leq} \left\Vert\Delta_{\hat{\mathbf{X}}_1}\mathbf{G}_{(1)}\left(\hat{\mathbf{X}}_1 \otimes \hat{\mathbf{X}}_2\right)^{\top}\right\Vert_F +  \left\Vert\Delta_{\hat{\mathbf{X}}_2}\mathbf{G}_{(2)}\left(\hat{\mathbf{X}}_1 \otimes{\mathbf{X}}_1\right)^{\top}\right\Vert_F\nonumber\\&\;\; + \left\Vert\Delta_{\hat{\mathbf{X}}_1}\mathbf{G}_{(1)}\left({\mathbf{X}}_2 \otimes {\mathbf{X}}_1\right)^{\top}\right\Vert_F\nonumber\\
    &\stackrel{(b)}{\leq} \left\Vert\Delta_{\hat{\mathbf{X}}_1}\mathbf{G}_{(1)}\right\Vert_F \left\Vert\hat{\mathbf{X}}_2\right\Vert\left\Vert\hat{\mathbf{X}}_1\right\Vert + \left\Vert\Delta_{\hat{\mathbf{X}}_2}\mathbf{G}_{(2)}\right\Vert_F \left\Vert{\mathbf{X}}_1\right\Vert\left\Vert\hat{\mathbf{X}}_1\right\Vert \nonumber\\&\;\;+ \left\Vert\Delta_{\hat{\mathbf{X}}_1}\mathbf{G}_{(1)}\right\Vert_F \left\Vert{\mathbf{X}}_2\right\Vert\left\Vert{\mathbf{X}}_1\right\Vert\nonumber\\
    &\stackrel{(c)}{\leq}(1+\tau)^2 \left\Vert\Delta_{\hat{\mathbf{X}}_1}\mathbf{\Sigma}_1\right\Vert_F + (1+\tau)  \left\Vert\Delta_{\hat{\mathbf{X}}_2}\mathbf{\Sigma}_2\right\Vert_F + \left\Vert\Delta_{\hat{\mathbf{X}}_1}\mathbf{\Sigma}_1\right\Vert_F\nonumber\\
    &\stackrel{(d)}{\leq} \gamma (1+\tau)^2 \mbox{dist}\left(\hat{\mathbf{F}}_{t},\mathbf{F}\right)\label{eq:g_1_1},
\end{align} \noindent {where $(a)$ and $(b)$ are due to the decomposition of $\mathcal{H}=\left(\mathbf{X}_1,\mathbf{X}_2,\mathbf{X}_1\right)\cdot\mathcal{G}$, the Tucker decomposition norm invariance $\| \mathcal{X}\|_F = \| \mathbf{X}_{(k)}\| _F \; \forall \; k$, and the Frobenius and spectral norm inequalities, together with $\|\mathbf A\mathbf B\|_F\leq
\|\mathbf A\|_F\|\mathbf B\|$,
$\|\mathbf A\otimes\mathbf B\|
=
\|\mathbf A\|\|\mathbf B\|$,
and $\|\mathbf X_k\|=1$. {Inequality $(c)$ holds since
$\mathbf G_{(k)}\mathbf G_{(k)}^\top=\mathbf\Sigma_k^2$, then
we have
$\|\Delta_{\hat{\mathbf X}_k}\mathbf G_{(k)}\|_F
=
\|\Delta_{\hat{\mathbf X}_k}\mathbf\Sigma_k\|_F$,
$\|\hat{\mathbf X}_k\|\leq 1+\tau$ from (\ref{eq:b1}) and the definition of $\tau$, $(d)$ follows from \eqref{eq:dist_def}.}}
Finally, plugging \eqref{eq:g_1_1} into \eqref{eq:g_1}, we have the bound for $\circledast_1$ as 
\begin{align}
    \circledast_1 \leq (1+\tau)^6 \left(\gamma^2(1+\tau)^4 \mbox{dist}^2\left(\hat{\mathbf{F}}_{t},\mathbf{F}\right)+ \left\Vert\mathbf{T}\right\Vert_{F}^2\right). \label{eq:circstar}
\end{align}  

By plugging  {\eqref{eq:circstar}} into \eqref{eq:hatg_g}, and combining 
the bounds for $\left\Vert\hat{\mathbf{X}}_{t+1,1}\mathbf{Q}_{t}-\mathbf{X}_{1}\right\Vert_{F}^{2}$, $\left\Vert\hat{\mathbf{X}}_{t+1,2}\mathbf{R}_{t}-\mathbf{X}_{2}\right\Vert_{F}^{2}$ and  $\left\Vert\left(\mathbf{Q}_{t}^{-1},\mathbf{R}_{t}^{-1},\mathbf{Q}_{t}^{-1}\right)\cdot\hat{\mathcal{G}}_{t+1}-\mathcal{G}\right\Vert_{F}^{2}$ in \eqref{eq:xq_x1} and \eqref{eq:hatg_g}, we can bound \eqref{eq:11} as 
\begin{align}
    &\mbox{dist}^{2}\left(\hat{\mathbf{F}}_{t+1},\mathbf{F}\right)\nonumber
    \leq  ((1-\eta)^2+\rho_0+\rho_1)\mbox{dist}^{2}\left(\hat{\mathbf{F}}_{t},\mathbf{F}\right)\\&+ 8\eta^2
(\epsilon_1K_1+qL_1)^2 + 8\eta^2
(\epsilon_1K_2+qL_2)^2\nonumber\\
    &+ \eta^2(1+\tau)^6 \left(\gamma^2(1+\tau)^4 \mbox{dist}^2\left(\hat{\mathbf{F}}_{t},\mathbf{F}\right)+ \left\Vert\mathbf{T}\right\Vert_{F}^2\right) \label{eq:final_bound}
\end{align} 
With $K=max\{K_1,K_2\}$ and $L=max\{L_1,L_2\}$,
\eqref{eq:final_bound} becomes
\begin{align}
    &\mbox{dist}^{2}\left(\hat{\mathbf{F}}_{t+1},\mathbf{F}\right)\\
    &\leq \left(1-2\eta+\eta^2+\rho_0+\rho_1+\eta^2\gamma^2(1+\tau)^{10}\right)\mbox{dist}^2\left(\hat{\mathbf{F}}_{t},\mathbf{F}\right)\nonumber\\& \;\;\;+ 16\eta^2(\epsilon_1K+qL)^2
+\eta^2(1+\tau)^6
\left\|
\mathbf{T}
\right\|_F^2\nonumber
\end{align} {Invoking the assumption on $\eta$ in the lemma statement,} $\frac{1-1/\sqrt{2}}{2\text{min}_{k}\sigma_{min}^2(\mathbf{\Sigma}_k)}< \eta< \text{max}\left\{\frac{1}{2\text{min}_{k}\sigma_{max}^2(\mathbf{\Sigma}_k)},\frac{1}{1+\gamma^2(1+\tau)^{10}}\right\}$, we have,
\begin{align}
    &\mbox{dist}^{2}\left(\hat{\mathbf{F}}_{t+1},\mathbf{F}\right) \nonumber
    \\&\leq \rho\cdot\mbox{dist}^2\left(\hat{\mathbf{F}}_{t},\mathbf{F}\right)+ \eta^2
\left[
16(\epsilon_1K+L\left\|
\mathbf{T}
\right\|_F^2)^2
+
(1+\tau)^6
\left\|
\mathbf{T}
\right\|_F^2
\right],\nonumber
\end{align} \noindent where $\rho<1$. \qed
\vspace{-0.15in}
\section{Proof of Lemma~\ref{lem:tentofac}}\label{prof:lem2}
\vspace*{-0.05in}
Given $\mathcal{H}=(\mathbf{X}_1,\mathbf{X}_2,\mathbf{X}_1)\cdot \mathcal{G}$ and a repeated addition of zero, we have, 
\begin{align}
&\left(\hat{{\mathbf{X}}}_{1},\hat{{\mathbf{X}}}_{2},\hat{\mathbf{X}}_{1}\right)\cdot{\hat{\mathcal{G}}}-\mathcal{H}= \left(\hat{{\mathbf{X}}}_{1},\hat{{\mathbf{X}}}_{2},\hat{\mathbf{X}}_{1}\right)\cdot{\hat{\mathcal{G}}} - \left(\hat{{\mathbf{X}}}_{1},\hat{{\mathbf{X}}}_{2},\hat{\mathbf{X}}_{1}\right)\cdot{{\mathcal{G}}}\nonumber \\&+ \left(\hat{{\mathbf{X}}}_{1},\hat{{\mathbf{X}}}_{2},\hat{\mathbf{X}}_{1}\right)\cdot{{\mathcal{G}}} - \left({{\mathbf{X}}}_{1},\hat{{\mathbf{X}}}_{2},\hat{\mathbf{X}}_{1}\right)\cdot{{\mathcal{G}}}+ \left({{\mathbf{X}}}_{1},\hat{{\mathbf{X}}}_{2},\hat{\mathbf{X}}_{1}\right)\cdot{{\mathcal{G}}}\nonumber \\&- \left({{\mathbf{X}}}_{1},{{\mathbf{X}}}_{2},\hat{\mathbf{X}}_{1}\right)\cdot{{\mathcal{G}}} + \left({{\mathbf{X}}}_{1},{{\mathbf{X}}}_{2},\hat{\mathbf{X}}_{1}\right)\cdot{{\mathcal{G}}} - \left({{\mathbf{X}}}_{1},{{\mathbf{X}}}_{2},{\mathbf{X}}_{1}\right)\cdot{{\mathcal{G}}}.\nonumber
\end{align}
Once again we let $\tau=\epsilon_1\sigma_{min}$. We consider the Frobenius norm,
{\begin{align}   &\left\Vert\left(\hat{{\mathbf{X}}}_{1},\hat{{\mathbf{X}}}_{2},\hat{\mathbf{X}}_{1}\right)\cdot{\hat{\mathcal{G}}}-\mathcal{H}\right\Vert_{F} \nonumber\\
&\stackrel{(a)}{\leq} \left\Vert\left(\hat{{\mathbf{X}}}_{1},\hat{{\mathbf{X}}}_{2},\hat{\mathbf{X}}_{1}\right)\cdot \Delta_{\hat{\mathcal{G}}}\right\Vert_{F} + \left\Vert\Delta_{\hat{\mathbf{X}}_1}\mathbf{G}_{(1)}\left(\hat{\mathbf{X}}_1\otimes \hat{\mathbf{X}}_2\right)^{\top}\right\Vert_{F}\nonumber\\
  &+ \left\Vert\Delta_{\hat{\mathbf{X}}_2}\mathbf{G}_{(2)}\left(\hat{\mathbf{X}}_1\otimes {\mathbf{X}}_1\right)^{\top}\right\Vert_{F} + \left\Vert\Delta_{\hat{\mathbf{X}}_1}\mathbf{G}_{(1)}\left({\mathbf{X}}_2\otimes {\mathbf{X}}_1\right)^{\top}\right\Vert_{F}\nonumber
          \end{align}  
\begin{align}  
    &\stackrel{(b)}{\leq}\left\Vert\hat{\mathbf{X}}_1\right\Vert\left\Vert\hat{\mathbf{X}}_2\right\Vert\left\Vert\hat{\mathbf{X}}_1\right\Vert\left\Vert\Delta_{\hat{\mathcal{G}}}\right\Vert_F+\left\Vert\Delta_{\hat{\mathbf{X}}_1}\mathbf{G}_{(1)}\right\Vert_F\left\Vert\hat{\mathbf{X}}_2\right\Vert\left\Vert\hat{\mathbf{X}}_1\right\Vert\nonumber\\
&+\left\Vert\Delta_{\hat{\mathbf{X}}_2}\mathbf{G}_{(2)}\right\Vert_F\left\Vert\hat{\mathbf{X}}_1\right\Vert\left\Vert{\mathbf{X}}_1\right\Vert + \left\Vert\Delta_{\hat{\mathbf{X}}_1}\mathbf{G}_{(1)}\right\Vert_F\left\Vert{\mathbf{X}}_2\right\Vert\left\Vert{\mathbf{X}}_1\right\Vert\nonumber\\
    &\stackrel{(c)}{\leq} (1+\tau)^3\left\Vert\Delta_{\hat{\mathcal{G}}}\right\Vert_F+\left(\tau^2+2\tau+2\right) \left\Vert\Delta_{\hat{\mathbf{X}}_1}\mathbf{\Sigma}_1\right\Vert_F+(1+\tau)\left\Vert\Delta_{\hat{\mathbf{X}}_2}\mathbf{\Sigma}_2\right\Vert_F\nonumber,\\
    &\stackrel{(d)}{\leq} \frac{4}{3}(1+\tau)^3\cdot\left(\left\Vert\Delta_{\hat{\mathbf{X}}_1}\right\Vert_F\left\Vert\mathbf{\Sigma}_1\right\Vert+\left\Vert\Delta_{\hat{\mathbf{X}}_2}\right\Vert_F\left\Vert\mathbf{\Sigma}_2\right\Vert+\left\Vert\Delta_{\hat{\mathcal{G}}}\right\Vert_F\right)\nonumber\\
    &\leq \frac{4}{3}(1+\tau)^3 \left(\left\Vert\Delta_{\hat{\mathbf{X}}_1}\right\Vert+\left\Vert\Delta_{\hat{\mathbf{X}}_2}\right\Vert+\left\Vert\Delta_{\hat{\mathcal{G}}}\right\Vert\right)\nonumber\\&\cdot\max\left\{\sigma_1\left(\mathbf{H}_{(1)}\right),\sigma_1\left(\mathbf{H}_{(2)}\right),1\right\}\stackrel{(e)}{\leq} \frac{4}{3}\left(1+\epsilon_1\sigma_{\min}\right)^3\zeta\cdot\text{dist}\left(\mathbf{F},\hat{\mathbf{F}}\right),\nonumber
\end{align} \noindent where $\zeta=\max\left\{\sigma_1\left(\mathbf{H}_{(1)}\right),\sigma_1\left(\mathbf{H}_{(2)}\right),1\right\}$. Inequality $(a)$ follows from the triangle inequality, {the mode-$k$ matricization identity of Tucker product defined in Section~\ref{sec:prelim}, for example, $\left((\Delta_{\hat{\mathbf{X}}_1},\hat{\mathbf{X}}_2,\hat{\mathbf{X}}_1)\cdot \mathcal{G}\right)_{(1)}=\Delta_{\hat{\mathbf{X}}_1}\mathbf{G}_{(1)}(\hat{\mathbf{X}}_1\otimes\hat{\mathbf{X}}_2)^{\top}$} and  the Tucker decomposition norm invariance $\| \mathcal{X}\|_F = \| \mathbf{X}_{(k)}\| _F \; \forall \; k$. Inequalities $(b)$ and $(c)$ use Lemma~\ref{lem:spec_bound} together with
$\|\hat{\mathbf X}_k\|
\le
1+\epsilon_1\sigma_{\min} = 1 + \tau$, $k=1,2$. Inequality $(d)$ follows from the submultuplicativity of the spectral and Frobenuis norms, {i.e., $\|\mathbf{AB}\|\leq \|\mathbf{A}\|_F\|\mathbf{B}\|_2$,} and $(e)$ derives from \eqref{eq:dist_def}.} 

\vspace{-0.17in}
\section{Proof of Lemma~\ref{lem:error}}
\vspace*{-0.05in}
\label{prof:lem3}
Note that $\left\Vert\mathbf{E}_{\qhats,(1)}\right\Vert_{F}^{2}$ is equivalent to $\left\Vert\mathcal{E}_{\qhats}\right\Vert_{F}^{2}$. Let us recall $\mathcal{E}_{\qhats} = \hat{\mathcal{Q}}\times_2\mathbf{S} - \mathcal{H}$, where $\hat{\mathcal{Q}}$ is obtained from \eqref{eq:Qhat}. Then we have, $\left\Vert\mathcal{E}_{\qhats}\right\Vert_{F}^2=\left\Vert\hat{\mathcal{Q}}\times_2\mathbf{S} - \mathcal{H}\right\Vert_{F}^2=\left\Vert\mathbf{H}_{(2)}-\mathbf{S}\hat{\mathbf{Q}}_{(2)}\right\Vert_{F}^2$.
The matrix $\hat{\mathbf{Q}}_{(2)}$ is obtained by solving the following problem, $\hat{\mathbf{Q}}_{(2)} = {\arg} \min_{\bar{\mathbf{Q}}_{(2)}}\;\left\Vert \mathbf{\Psi}\mathbf{H}_{(2)}-\mathbf{\Psi}\mathbf{S}\bar{\mathbf{Q}}_{(2)}\right\Vert_{F}^{2}$.
Denoting the solution to the least squares problem $\min_{\Tilde{\mathbf{Q}}_{(2)}}\;\left\Vert \mathbf{H}_{(2)}-\mathbf{S}\Tilde{\mathbf{Q}}_{(2)}\right\Vert_{F}^{2}$  by $\mathbf{Q}_{(2)}^{*}$, { we have that $\mathbf{H}_{(2)} = \mathbf{S} \mathbf{Q}_{(2)}^{*}+\mathbf{H}_{(2)}^{\perp}$.}  We next plug this value into the desired sketched regression problem,
\vspace*{-0.05in}
\begin{align}
    &\min_{\bar{\mathbf{Q}}_{(2)}}\;\left\Vert \mathbf{\Psi}\mathbf{H}_{(2)}-\mathbf{\Psi}\mathbf{S}\bar{\mathbf{Q}}_{(2)}\right\Vert_{F}^{2}\nonumber\\&= \min_{\bar{\mathbf{Q}}_{(2)}}\;\left\Vert \mathbf{\Psi}\left(\mathbf{S}\mathbf{Q}_{(2)}^{*}+\mathbf{H}_{(2)}^{\perp}\right)-\mathbf{\Psi}\mathbf{S}\bar{\mathbf{Q}}_{(2)}\right\Vert_{F}^{2}\nonumber\\
    &\stackrel{(a)}{=}\min\;\left\Vert\mathbf{\Psi}\mathbf{H}_{(2)}^{\perp}+\mathbf{\Psi}\mathbf{U}_{S}\left(\mathbf{Y}_{(2)}^{*}-\bar{\mathbf{Y}}_{(2)}\right)\right\Vert_{F}^{2},\nonumber
\end{align}\noindent {where $(a)$ follows from $\mathbf{S}\mathbf{Q}_{(2)}^{*}=\mathbf{U}_{S}\mathbf{Y}_{(2)}^{*}$, and the same relationship holds for $\bar{\mathbf{Q}}_{(2)}$, $\bar{\mathbf{Y}}_{(2)}$, i.e., $\mathbf{S}\bar{\mathbf{Q}}_{(2)}=\mathbf{U}_{S}\bar{\mathbf{Y}}_{(2)}$ and $\hat{\mathbf{Q}}_{(2)}$,$\hat{\mathbf{Y}}_{(2)}$, i.e., $\mathbf{S}\hat{\mathbf{Q}}_{(2)}=\mathbf{U}_{S}\hat{\mathbf{Y}}_{(2)}$.} $\hat{\mathbf{Y}}_{(2)}$ {solves the following least squares problem}:
\vspace*{-0.05in}
{\begin{align}
    \hat{\mathbf Y}_{(2)}=
\arg\min_{\bar{\mathbf Y}_{(2)}}
\left\|
\mathbf\Psi\mathbf H_{(2)}
-
\mathbf\Psi\mathbf U_S\bar{\mathbf Y}_{(2)}
\right\|_F^2.
\end{align}
The first-order optimality condition yields the following normal equation $\left(\mathbf\Psi\mathbf U_S\right)^\top
\left(\mathbf\Psi\mathbf U_S\hat{\mathbf Y}_{(2)}
-\mathbf\Psi\mathbf H_{(2)}\right)=
\mathbf 0.$
Using $\mathbf{H}_{(2)}
=
\mathbf {U}_S\mathbf Y_{(2)}^*
+
\mathbf {H}_{(2)}^\perp,$
we obtain \vspace{-0.05in}
\begin{align}
&
\left(\mathbf\Psi\mathbf U_S\right)^\top
\mathbf\Psi\mathbf U_S
\left(
\hat{\mathbf Y}_{(2)}
-
\mathbf Y_{(2)}^*
\right)=
\left(\mathbf\Psi\mathbf U_S\right)^\top
\mathbf\Psi\mathbf H_{(2)}^\perp.
\end{align}}
\noindent 
{Assuming 
$\sigma_{\min}^{2}\left(\mathbf{\Psi}\mathbf{U}_{S}\right) \geq 1/\sqrt{2}$,} we have that $\sigma_{i}\left(\left(\mathbf{\Psi}\mathbf{U}_{S}\right)^{\top}\mathbf{\Psi}\mathbf{U}_{S}\right)=\sigma_{i}^{2}\left(\mathbf{\Psi}\mathbf{U}_{S}\right)\geq 1/\sqrt{2}$. Thus, we have
\begin{align}
    \left\Vert\hat{\mathbf{Y}}_{(2)}-\mathbf{Y}_{(2)}^{*}\right\Vert^2_{F} \slash 2 &\leq \left\Vert\left(\mathbf{\Psi}\mathbf{U}_{S}\right)^{\top}\mathbf{\Psi}\mathbf{U}_{S}\left(\hat{\mathbf{Y}}_{(2)}-\mathbf{Y}_{(2)}^{*}\right)\right\Vert_{F}^{2}\nonumber
        \end{align}
    \begin{align}& \stackrel{(a)}{=} \left\Vert\left(\mathbf{\Psi}\mathbf{U}_{S}\right)^{\top}\mathbf{\Psi}\mathbf{H}_{(2)}^{\perp}\right\Vert_{F}^{2},\nonumber
\end{align} \noindent where $(a)$ is due to $\mathbf{U}_{S}\left(\hat{\mathbf{Y}}_{(2)}-\mathbf{Y}_{(2)}^{*}\right) = \mathbf{H}_{(2)}^{\perp}$. Finally, we apply Lemma~\ref{lem:mat_approx} to the right hand side of the inequality to obtain:
\begin{align}\label{eq:yhat_y}
    \left\Vert\hat{\mathbf{Y}}_{(2)}-\mathbf{Y}_{(2)}^{*}\right\Vert^2_{F} \slash 2 \leq \epsilon_2\left\Vert\mathbf{H}_{(2)}^{\perp}\right\Vert_{F}^{2}.
\end{align} This result implies:
\begin{align}
    &\left\Vert\mathbf{H}_{(2)}-\mathbf{S}\hat{{\mathbf{Q}}}_{(2)}\right\Vert_{F}^2= \left\Vert\mathbf{H}_{(2)}- \mathbf{S}\mathbf{Q}_{(2)}^{*}+\mathbf{S}\mathbf{Q}_{(2)}^{*}- \mathbf{S}\hat{{\mathbf{Q}}}_{(2)}\right\Vert_{F}^2\nonumber\\
    &\leq\left\Vert\mathbf{H}_{(2)}- \mathbf{S}\mathbf{Q}_{(2)}^{*}\right\Vert_{F}^{2}+\left\Vert\mathbf{S}\mathbf{Q}_{(2)}^{*}- \mathbf{S}\hat{{\mathbf{Q}}}_{(2)}\right\Vert_{F}^{2}\nonumber\\
    &\stackrel{(a)}{=}\left\Vert\mathbf{H}_{(2)}^{\perp}\right\Vert_{F}^{2} + \left\Vert\mathbf{U}_{S}\left(\mathbf{Y}_{(2)}^{*}-\hat{\mathbf{Y}}_{2}\right)\right\Vert_{F}^2\nonumber\\
    &\stackrel{(b)}{\leq} \left\Vert\mathbf{H}_{(2)}^{\perp}\right\Vert_{F}^{2} + 2\epsilon_2\left\Vert\mathbf{H}_{(2)}^{\perp}\right\Vert_{F}^{2}= \left(1+2\epsilon_2\right)\left\Vert\mathbf{H}_{(2)}^{\perp}\right\Vert_{F}^{2}\nonumber \\&= \left(1+2\epsilon_2\right)\left\Vert\mathbf{H}_{(2)}-\mathbf{S}\mathbf{Q}_{(2)}^{*}\right\Vert_{F}^{2},\label{eq:H_SQ}
\end{align} \noindent {where $(a)$ is from $\mathbf{S}\mathbf{Q}_{(2)}^{*}=\mathbf{U}_{S}\mathbf{Y}_{(2)}^{*}$ and $\mathbf{S}\hat{\mathbf{Q}}_{(2)}=\mathbf{U}_{S}\hat{\mathbf{Y}}_{(2)}$ and  $(b)$ is due to \eqref{eq:yhat_y}. } Next, we further bound $\left\Vert\mathbf{H}_{(2)}-\mathbf{S}\mathbf{Q}_{(2)}^{*}\right\Vert_{F}^{2}$ as $\left\Vert\mathbf{H}_{(2)}-\mathbf{S}\mathbf{Q}_{(2)}^{*}\right\Vert_{F}^{2} = \left\Vert\mathbf{H}_{(2)}^{\perp}\right\Vert_{F}^{2}= \left\Vert\mathbf{H}_{(2)}\right\Vert_{F}^{2}-\left\Vert\mathbf{U}_{S}\mathbf{U}_{S}^{\trasp}\mathbf{H}_{(2)}\right\Vert_{F}^{2}$. 
{Recall $\kappa \doteq
\frac{
\left\|
\mathbf U_S\mathbf U_S^\top\mathbf H_{(2)}
\right\|_F
}{
\left\|
\mathbf H_{(2)}
\right\|_F
},\;\;
\kappa\in(0,1]$. We have
\begin{align}
&\left\Vert\mathbf{H}_{(2)}\right\Vert_{F}^{2}-\left\Vert\mathbf{U}_{S}\mathbf{U}_{S}^{\trasp}\mathbf{H}_{(2)}\right\Vert_{F}^{2}\nonumber\\
    &= \kappa^{-2}\left\Vert\mathbf{U}_{S}\mathbf{U}_{S}^{\trasp}\mathbf{H}_{(2)}\right\Vert_{F}^{2} - \left\Vert\mathbf{U}_{S}\mathbf{U}_{S}^{\trasp}\mathbf{H}_{(2)}\right\Vert_{F}^{2}\nonumber\\
    &=\left(\kappa^{-2}-1\right)\left\Vert\mathbf{S}\mathbf{Q}_{(2)}^{*}\right\Vert_{F}^{2}\stackrel{(a)}{\leq} \left(\kappa^{-2}-1\right)\sigma_{max}^{2}\left(\mathbf{Q}_{(2)}^{*}\right)\left\Vert\mathbf{S}\right\Vert_{F}^{2}, \nonumber
\end{align}} \noindent Inequality $(a)$ follows from the norm inequality.  By plugging the above equation in \eqref{eq:H_SQ}, we complete the proof. \qed

\vspace{-0.1in}


\bibliography{ref,ref_tensor}

\end{document}